\documentclass[11pt]{amsart}

\usepackage{amsmath, amssymb, amsfonts, amsthm}
\usepackage{graphicx}
\usepackage{xcolor}
\usepackage{bm}
\usepackage{subfigure}
\usepackage{graphicx}
\usepackage{mathabx}
\usepackage{multirow}
\usepackage{setspace}
\usepackage[title]{appendix}
\usepackage{mathrsfs} 
\usepackage{enumerate}
\usepackage{arydshln}
\usepackage{mathdots}
\usepackage{upgreek}
\usepackage{comment}
\usepackage{soul}

\setstcolor{magenta}

\usepackage{amsthm}
\theoremstyle{plain}
\newtheorem{theorem}{Theorem}[section]

\theoremstyle{plain}
\newtheorem{lemma}[theorem]{Lemma}

\theoremstyle{plain}
\newtheorem{Prop}[theorem]{Proposition}

\theoremstyle{plain}
\newtheorem{corollary}[theorem]{Corollary}

\theoremstyle{plain}

\theoremstyle{definition}
\newtheorem{Def}[theorem]{Definition}
\theoremstyle{remark} 
\newtheorem{remark}[theorem]{Remark}

\theoremstyle{definition}

\theoremstyle{definition} 
\newtheorem{Example}[theorem]{Example}

\usepackage{tikz}
\usepackage{tikz-cd}
\usepackage{wrapfig}
\usepackage{float}

\usetikzlibrary{matrix,arrows,decorations.pathmorphing}

\def\ed{\mathrm{d}}

\def\I{\mathcal I}

\def\R{\mathbb R}

\def\x{{\bf x}}
\def\f{{\bf f}}
\def\u{{\bf u}}
\def\v{{\bf v}}

\def\y{{\bf y}}

\def\ut{\uptau}

\def\rank{{\rm rank } }

\def\alg{{\rm alg}}

\def\W{\wedge}
\def\<{\langle}
\def\>{\rangle}

\def\pr{{\rm pr}}

\def\({\left(}
\def\){\right)}

\def\lbb{[\![}
\def\rbb{]\!]}

\def\bs{\boldsymbol}

\DeclareFontFamily{U}{MnSymbolC}{}
\DeclareSymbolFont{MnSyC}{U}{MnSymbolC}{m}{n}
\DeclareFontShape{U}{MnSymbolC}{m}{n}{
    <-6>  MnSymbolC5
   <6-7>  MnSymbolC6
   <7-8>  MnSymbolC7
   <8-9>  MnSymbolC8
   <9-10> MnSymbolC9
  <10-12> MnSymbolC10
  <12->   MnSymbolC12}{}
\DeclareMathSymbol{\lefthook}{\mathbin}{MnSyC}{'270}

\tikzset{
  symbol/.style={
    draw=none,
    every to/.append style={
      edge node={node [sloped, allow upside down, auto=false]{$#1$}}}
  }
}

\newcommand{\eqdef}{\overset{\mathrm{}}{=\joinrel=}}

\title[On Absolute Equivalence and Linearization I]{On Absolute Equivalence and Linearization I}

\author{Jeanne N. Clelland}
\address{Department of Mathematics, 395 UCB, University of
Colorado, Boulder, CO 80309-0395}
\email{Jeanne.Clelland@colorado.edu}

\author{Yuhao Hu}
\address{Key Laboratory of Pure and Applied Mathematics, 
School of Mathematical Sciences, Peking University, Beijing, 100871, P. R. China
}
\email{yuhao.hu@math.pku.edu.cn, yuhao.hu@colorado.edu}



\thanks{The first author was supported by a Collaboration Grant for Mathematicians from the Simons Foundation, and 
the second author was partially supported by the China Postdoctoral Science Foundation Grant 2021TQ0014.}

\begin{document}
\maketitle

\begin{abstract}
	In this paper, we study the absolute equivalence between Pfaffian systems
	with a degree $1$ independence condition and obtain 
	 structural results, particularly for systems of corank $3$. We apply these 
	 results to understanding dynamic feedback linearization of control systems with $2$ inputs.
\end{abstract}

%\setcounter{tocdepth}{1}
%\tableofcontents

%Intro
\section{Introduction}

The notion of \emph{absolute equivalence} between two differential systems was introduced by \'Elie Cartan in
\cite{Cartan14}. 

By the time of Cartan's writing, Hilbert \cite{Hilbert35} asked the question: \emph{When can the
general solutions of an ODE for two unknown functions be expressed in a determined way 
in terms of an arbitrary function and its successive derivatives?}
He proved that such a property is not enjoyed by the ODE
		\begin{equation}\label{HilbertCartanEq}
			\frac{\ed x}{\ed t} = \(\frac{\ed^2y}{\ed t^2}\)^2.
		\end{equation}
This is in contrast with, for example, the ODE
		\[
			\(\frac{\ed x}{\ed t}\)^2 + \(\frac{\ed y}{\ed t}\)^2 = 1,
		\]
whose general (real analytic) solutions can be expressed as
		\begin{align*}
			t& = f''(\alpha) + f(\alpha),\\
			\(\begin{array}{c}
				x\\
				y
			\end{array}	
			\)	& = \(\begin{array}{cc}
				\sin \alpha & \cos\alpha\\
				-\cos\alpha & \sin\alpha
			\end{array}\)\(\begin{array}{c}
				f'(\alpha)\\
				f''(\alpha)
			\end{array}\)
		\end{align*}
for an arbitrary function $f(\alpha)$. 

Cartan realized that a key to answering Hilbert's question was to understand when one
can establish a one-to-one correspondence between the solutions of two differential systems
$\mathcal{E}_1$ and $\mathcal{E}_2$; and he called $\mathcal{E}_1$ and $\mathcal{E}_2$
``absolutely equivalent" if a particular kind of such correspondence exists. 

Furthermore, the differential systems that enter Hilbert's question are special cases of \emph{rank-$n$ Pfaffian systems in $n+2$ variables with a degree $1$ independence condition}\footnote{We will call
such a Pfaffian system a \emph{system of corank $2$}.}.
Given such a Pfaffian system $I$, one can compute its successive \emph{derived systems} $I^{(k)}$. 
The integers $r_k:=\rank\(I^{(k)}/I^{(k+1)}\)$ are important invariants of $I$ under diffeomorphisms. Only two cases can occur:
either $r_k\le 1$ for all $k$, in which case $I$ is said to have \emph{class $0$}; or $r_\ell = 2$ for the first time for some $\ell$ ($\ell\ge 1$), in which case $I$ is said to have \emph{class $\rank(I^{(\ell)})$}. (See \cite{gardner67} for details.)
In the latter case, Cartan called $I^{(\ell-1)}$
the \emph{normal system} of $I$.  Furthermore, he proved the following theorem.
\begin{theorem} \cite{Cartan14}
Two systems $I$, $\bar I$ of corank $2$ are absolutely equivalent if and only if the following two conditions hold
\begin{enumerate}[\bf i.]
	\item{$I$ and $\bar I$  have the same class;}
	\item{when both systems have class $0$, 
	the terminal derived systems $I^{(\infty)}$, $\bar I^{(\infty)}$ satisfy $\rank(I^{(\infty)}) = \rank(\bar I^{(\infty)}) $;
		 when both systems have the same positive class, the corresponding normal systems 
		$I^{(\ell-1)}$ and $\bar I^{(\bar \ell - 1)}$ can be transformed into each other
		 by a diffeomorphism between their underlying manifolds.}
\end{enumerate}
\end{theorem}
In terms of these, $I$ satisfies Hilbert's condition precisely when $I$ is absolutely equivalent to the Pfaffian system
generated by a single $1$-form
	\[
		\ed x - u\ed t,
	\]
which holds if and only if $I$ has class $0$ and $I^{(\infty)} = 0$. That the Pfaffian system corresponding to 
\eqref{HilbertCartanEq} does not satisfy this condition can be checked simply by computing 
derived systems.

Cartan's motivation was apparently classical differential geometry, 
as the second half of \cite{Cartan14} shows.
This serves as an interesting contrast, as eighty years later,
a significant amount of interest returned to this paper of Cartan, this time
motivated by applications in control theory, represented by the works
\cite{shadwick1990absolute}, \cite{gardner1990feedback} and \cite{Sluis}, to mention a few.

An \emph{autonomous control system} in $n$ states and $m$ inputs (aka. controls) $(n\ge m)$ is an ODE system of the form
	\[
		\dot\x = \f(\x,\u),\qquad \x\in \Omega\subseteq\R^n,~ \u\in \R^m.
	\]
In addition, we will assume that
	\[
		\rank\(\frac{\partial f^i}{\partial u^\alpha}\) = m.
	\]	
The system is called {\em controllable} if, given any two states $\x_1, \x_2 \in \R^n$, there exists a solution $\x(t)$ and $t_1, t_2 \in \R$ such that $\x(t_1) = \x_1$ and $\x(t_2) = \x_2$.  

The simplest class of such systems are those that are \emph{linear}:
	\[
		\dot\x = A\x + B\u,
	\]
where $A, B$ are constant matrices with $B$ attaining full rank. Controllability in this case is characterized by \emph{Kalman's maximal rank condition}:
	\[
		\rank\(B|AB|A^2B|\cdots|A^{n-1}B\) = n.
	\]

It is well-known that a controllable linear system can be put into a {\em Brunovsk\'y normal form} \cite{Brunovsky}
by a linear \emph{feedback transformation} of the form
	\[
		\y = C\x,\qquad \v = D\x + E\u,
	\] 
where $C, D, E$ are constant matrices.
This normal form is defined as follows:
\begin{Def}
	A control system of the following form
\begin{equation}\label{Brunovsky-normal-form}
	\begin{alignedat}{1}
		\dot x^1_j &= x^1_{j+1}\quad (j = 0,\ldots, r_1),\\
		\dot x^2_j &= x^2_{j+1}\quad (j = 0,\ldots, r_2),\\
		&\cdots\\
		\dot x^p_j & = x^p_{j+1}\quad (j = 0, \ldots, r_p),
	\end{alignedat}	
\end{equation}
	where $r_1,\ldots, r_p\ge 0$, is called a \emph{Brunovsk\'y normal form}. In particular, $x^i_j$ 
	$(i = 1,\ldots, p; j = 0,\ldots, r_i)$ are the state variables; $x^i_{r_i+1}$ $(i = 1,\ldots, p)$ are the control variables.  
\end{Def}	
The Brunovsk\'y system \eqref{Brunovsky-normal-form} is particularly simple, in that its general solution may be expressed in terms of $p$ arbitrary functions $f^1(t), \ldots, f^p(t)$ as follows:
\begin{equation*}
\begin{alignedat}{3}
x^1_0(t) & = f^1(t),& \quad x^1_1(t) & = (f^1)'(t), & \quad \ldots \quad x^1_{r_1}(t) & = (f^1)^{(r_1)}(t), \\
&~\vdots & &~\vdots& & \\
x^p_0(t) & = f^p(t),& \quad x^p_1(t) & = (f^p)'(t), & \quad \ldots \quad x^1_{r_p}(t) & = (f^p)^{(r_p)}(t).
\end{alignedat}
\end{equation*}

More generally, it is also known when an arbitrary controllable autonomous system can be put
into a Brunovsk\'y normal form by a $t$-independent, not necessarily linear, feedback transformation:
	\[
		\y = \phi(\x),\qquad \v = \psi(\x,\u),
	\]
where $\y,\v$ are the new states and inputs, respectively. In fact,
 in \cite{gardner1992gs} (see also \cite[p.78, Theorem 35]{Sluis}), Gardner and Shadwick proved:
 \begin{quote}
  \emph{An autonomous control system, formulated
 as a Pfaffian system $I$ with the independence condition $\tau = \ed t$, 
 is locally feedback equivalent to a Brunovsk\'y normal form if and only if
 	\begin{enumerate}[\bf i.]
		\item{the terminal derived system $I^{(\infty)}$ vanishes;}
		\item{for each integer $k\ge 0$, the Pfaffian system generated by the $k$-th derived system $I^{(k)}$ 
		and $\ed t$ is Frobenius.}
	\end{enumerate}}	
\end{quote}

A more general question is: \emph{Given two control systems, how to tell whether they are absolutely equivalent?}
The answer to this question largely depends on the number of control variables that one is dealing with.
For control systems with a single control, the question reduces to the case addressed by \cite{Cartan14}.
A systematic study of cases with more than one control was carried out in the 1994 thesis \cite{Sluis} of Sluis.
In that thesis, Sluis followed Cartan, realizing that the key to understanding absolute equivalence
is to understand the so-called ``Cartan prolongations". A notable theorem that Sluis proved is
Theorem \ref{SluisThm} below, which allows one to understand any Cartan prolongation in terms of ``prolongations by differentiation". It is natural to wonder whether, in certain cases, the theorem can be used to classify 
control systems that are absolutely equivalent to a Brunovsk\'y normal form via a succession of such prolongations. This is a main motivation for the current work:  In this paper we consider systems with $n$ states and 2 controls, and our main result (Theorem \ref{typen2LinTheorem}) gives necessary and sufficient conditions for such a system to be dynamic feedback linearizable after a particular, {\em fixed} number $K$ of differentiations.  
Additionally, in Theorem \ref{typen2StrEqn} we derive structure equations that provide a starting point for {\em classifying} those control systems that can be linearized via a succession of $K$ differentiations.

In future work, we hope to prove an upper bound on the maximum number of differentiations required in order to perform such a linearization. This would then provide a complete list of necessary and sufficient conditions for such a system to be dynamic feedback linearizable via a succession of prolongations by differentiation.  Ideally, the structure equations of Theorem \ref{typen2StrEqn} might then be used to produce a complete classification of such systems.

This paper is organized as follows.

To start with (Sections \ref{cartanprolSec} and \ref{dynLinSec}), we remind the reader of various notions of prolongation of a system (Section \ref{ProlongationSection}) and the extension theorem of Sluis (Section \ref{extThmsec}). 
Then, we introduce the notion of a \emph{relative extension} (Section \ref{relativeExtsec}), which is 
a canonical construction that one can obtain from a Cartan prolongation. 
A relative extension can be viewed as a nested array of Pfaffian bundles.
In general, it need not induce a succession of Cartan prolongations. However, 
when the original Cartan prolongation is \emph{regular} and when the systems involved have corank $3$, the relative 
extensions do induce a succession of Cartan prolongations (Theorem \ref{corank3Thm}). This allows us
to relate various notions of equivalence obtained from different notions of prolongation (Theorem \ref{RrelAbsEq}),
which in turn allows us to describe in a canonical way the necessary and sufficient conditions for a 
type $(n,2)$ (that is, $n$ states and $2$ controls) control-type system to be dynamic feedback linearizable via a prolongation by differentiation involving $K$ differentiations
(Theorem \ref{typen2LinTheorem}). We then set up the structure equations (Theorem \ref{typen2StrEqn}) 
that are associated to a linearization. In a sequel to this paper, these structure equations will serve as 
the basis of classifying control systems that can be linearized after a particular, fixed number of prolongations.

%Symbols
\subsection{Symbols and Abbreviations}
\begin{quote}
{\begin{itemize}
	\setlength\parskip{0.5mm}
	\item[EDS:]{Exterior differential system(s).}

	\item[CTS:]{Control-type system(s).}
	
	\item[$(M,\I)$:]{An EDS with manifold $M$ and differential ideal $\I\subset \Omega^*(M)$.}
	
	\item[$\<\eta^1,\ldots,\eta^k\>$:]{The differential ideal in $\Omega^*(M)$ generated by differential forms $\eta^1,\ldots,\eta^k\in \Omega^*(M)$ and their exterior derivatives.}
							
	\item[$\<\eta^1,\ldots,\eta^k\>_\alg$:]{The ideal in $\Omega^*(M)$ algebraically generated 
							by differential forms $\eta^1,\ldots,\eta^k\in \Omega^*(M)$.}

	\item[$(\pr^{(k)}M, \pr^{(k)}\I)$:] {The $k$-th (total) prolongation of an EDS $(M,\I)$.}
	
	\item[$(M,I)$:]{A Pfaffian system determined by a vector subbundle $I\subset T^*M$.}
	
	\item[$(\pr^{(k)}M, \pr^{(k)}I)$:] {The $k$-th (total) prolongation of a Pfaffian system $(M,I)$.}
	
	\item[$I^{(k)}$:]{The $k$-th derived system of a Pfaffian system $I$.}

	\item[$\mathcal{C}(\I)$ (resp. $\mathcal{C}(I)$):] {The Cartan system of a differential ideal $\I\subset \Omega^*(M)$ (resp., of a Pfaffian system $I\subset T^*M$). By definition, it is the Frobenius system defined
	on $M$ whose
			leaves are precisely the Cauchy characteristics of $(M,\I)$ (resp. $(M,I)$).
		}
	
	\item[$\mathcal{C}(\omega^1,\ldots,\omega^k)_{\alg}$:] {The algebraic Cartan system associated
					to an algebraic ideal of $\Omega^*(M)$ generated by differential forms $\omega^1,\ldots, \omega^k$.
					By definition, this is the Pfaffian system dual to the distribution spanned by 
					all vector fields $X$ that satisfy 
					$X\lefthook \omega^i\in \<\omega^1,\ldots,\omega^k\>_\alg$
					for all $i= 1,\ldots,k$.
					}
	\item[$\bs\theta$:]{A vector valued $1$-form $(\theta^1,\ldots,\theta^k)^T$.}							
	\item[$\lbb\theta^1,\ldots,\theta^k\rbb$:]{The rank $k$ subbundle of $T^*M$, or the 
					corresponding Pfaffian system, spanned by
					$k$ linearly independent $1$-forms $\theta^1,\ldots,\theta^k$.}		
	\item[$\lbb I,\eta^1,\ldots,\eta^k\rbb$:]{The subbundle of $T^*M$ spanned by the sections of a subbundle $I\subset T^*M$ and differential $1$-forms $\eta^1,\ldots,\eta^k\in\Omega^1(M)$. Sometimes we write 
	 $\lbb I, \eta^i\rbb_{i = 1}^k$ for brevity.}											
	\item[$I_k$:] {The $k$-th extension of a system $I$ relative to a Cartan prolongation.}
	\item[$|\mathcal{S}|$:]{The cardinality of a finite index set $\mathcal{S}$.}		
		
\end{itemize}}
\end{quote}

%Cartan Prolongations
\section{The Structure of a Cartan Prolongation}\label{cartanprolSec}
	
	\begin{Def}\label{systemDef}
		A \emph{system} is a triple $(M,I;\tau)$, where $(M,I)$ is a Pfaffian system (i.e., $M$ is a smooth manifold, and  $I\subset T^*M$
		is a vector subbundle), and $\tau$, called the \emph{independence condition}, is an exact $1$-form which nowhere belongs to $I$, such that
	\begin{enumerate}[\bf i.]
		\item{$(M,I)$ admits no Cauchy characteristics, and}
		\item{$(M,I)$ admits no integral surface $\iota: S\hookrightarrow M$ that satisfies $\iota^*\tau \ne 0$.}
	\end{enumerate}
	\end{Def}
	
	This definition is local. 
	Generally, one would not require $\tau$ to be 
	exact (see \cite[p.35]{Sluis}), but an exact independence condition can always be chosen by shrinking 
	$M$ and allowing a negligible change in
	the set of integral curves of the Pfaffian system. 
	Having an exact independence condition distinguished will
	be convenient when we later work with control systems, for which $\ed t$, the differential of 
	time, is the natural choice of an independence condition.
	
	Furthermore, we will always assume that a system and its derived systems have constant ranks. This
	can be achieved
	by shrinking $M$, if needed. 
	
	%%Comment
	%{\color{red}\bf  2. The following sentence is new.} {\color{blue} (Comment: slightly edited)}
	
	Our motivation for including condition {\bf ii} in Definition \ref{systemDef} will be explained in
	Remark \ref{cauchyCharRemark} (Section \ref{ProlongationSection}).
	
	\begin{Def}
		A system $(M,I;\tau)$ is said to have \emph{type} $(n,m)$ if
		\begin{equation*}\label{typeDef}\left\{
		\begin{alignedat}{1}
			\dim(M)& = n+m+1,\\
			 \rank(I) &= n.
		\end{alignedat}\right. 
		\end{equation*}	 
			  The integer $m+1$ will be called the \emph{corank} of $(M,I;\tau)$.
	\end{Def}
	
	\begin{Def}
	Given two systems $(M,I;\tau)$ and $(\bar M,\bar I;\bar \tau)$, we say that they are
	 \emph{$\ut$-equivalent}
	if there exists a diffeomorphism
	\[
		\phi: M \rightarrow \bar M
	\]
	such that
	\[
		\phi^*\bar I = I, \qquad \phi^*\bar \tau  = \tau.
	\]
	\end{Def}

	{\remark Replacing the condition $\phi^*\bar\tau = \tau$ by
	 \[\phi^*\bar\tau \equiv \rho\tau \mod I\] 
	  for some nonzero function $\rho$,
	 one would obtain a more general notion of equivalence (\cite[p.38]{Sluis}).
	 Because of this, in the definition above, ``$\uptau$-" is used in order to eliminate possible confusion.}\\
	
	Sluis proved the following result.
	
	\begin{theorem} \cite[p.37, Theorem 14]{Sluis}\label{SluisCTSTheorem}
		A system $(M,I;\tau)$ admits local coordinates in which it corresponds to a (time-varying) control system
			\[
				\dot \x  = \f(\x,\u,t)
			\]
		such that $\tau = \ed t$ if and only if the Pfaffian system $\lbb I,\tau\rbb$ is Frobenius.
	\end{theorem}
	
	Motivated by this, we make the definition below.
	
	\begin{Def}
		A system $(M,I;\tau)$ of type $(n,m)$ is said to be a \emph{control-type system} (CTS) if $n \geq m$ and
		the Pfaffian system $\lbb I,\tau\rbb$ is Frobenius.
	\end{Def}
	
	{\remark \begin{enumerate}[\bf A.]
		\item{A type $(n,m)$ CTS corresponds to a control system with $n$ states and $m$ inputs. 
		}
		\item{Given two CTS $(M,I;\tau)$ and $(\bar M,\bar I;\bar\tau)$, 
		represented in coordinates by $\dot \x = \f(\x,\u,t)$ and $\dot {\bar\x}= \bar \f(\bar \x,\bar\u,\bar t)$, respectively, it is not difficult to see that they are $\ut$-equivalent as systems if and only if there exists
		a diffeomorphism of the form
		\begin{equation}\label{tpresTrans}
			\phi: (\x,\u,t)\mapsto (\bar\x,\bar\u, \bar t) = \(\bs\xi(\x,t),\bs\eta(\x,\u,t), t + T\),
		\end{equation}
		where $T$ is a constant, that transforms the second system into the first, and vice versa.	}
		\end{enumerate}
	}
	
	\subsection{Prolongations of a System}\label{ProlongationSection}
	Given a system $(M,I;\tau)$, one can define the following four types of prolongations: \emph{total prolongation},  \emph{partial prolongation}, \emph{prolongation by differentiation} and \emph{Cartan prolongation}.
	
	\begin{Def}\label{ProlongationsDef}
		Let $(M,I;\tau)$ be a system. As it will be sufficient for the this paper,
		 the following notions of prolongation are
		defined locally.
	\begin{enumerate}[\bf 1.]
	\item{
	By shrinking $M$, if needed, choose coordinates $(x^1,\ldots,x^n,u^1,\ldots,u^m,t)$ on $M$
	such that $I, \ed u^\alpha, \tau$ span the entire cotangent bundle $T^*M$.
	Let $\pr^{(1)}M:= M\times \R^m$, with $(\lambda^\alpha)$ being coordinates on the $\R^m$-component;
	  let
	\[	
		\pr^{(1)}I:= \lbb I, \ed u^\alpha - \lambda^\alpha \tau\rbb_{\alpha  = 1}^m,
	\]
	and let the pull-back of $\tau$ to $\pr^{(1)}M$ be denoted by the same letter.
	The system $(\pr^{(1)}M, \pr^{(1)}I;\tau)$ is called the \emph{total prolongation} of $(M,I;\tau)$. (Note that this system is usually just called the {\em prolongation} of $(M,I;\tau)$; we have added the adjective ``total" to distinguish it from the other types of prolongations being considered.)
	This definition is independent of the choice of local coordinates.}
	
	\item{A non-canonical way to ``prolong" $(M,I;\tau)$ is 
	by choosing $\mu^\alpha\in \Omega^1(M)$, $\alpha\in \mathcal{S}\subseteq \{1,2,\ldots,m\}$,
	linearly independent modulo $I$ and $\tau$, then letting
	$N:=M\times \R^{|\mathcal{S}|}$ and $J:=\lbb I,\mu^\alpha - \lambda^\alpha \tau\rbb_{\alpha\in \mathcal{S}}$, 
	where $(\lambda^\alpha)$ are 
	coordinates on the $\R^{|\mathcal{S}|}$-factor. The system $(N,J;\tau)$ is called a \emph{partial prolongation} 
	of $(M,I;\tau)$.
	}
	
	\item{Let $(x^1,\ldots,x^n, u^1,\ldots,u^m,t)$ be coordinates on $M$ such that $\ed t = \tau$ and that 
	$I, \ed u^\alpha, \ed t$ generate
	the entire $T^*M$. The system
		\[
			(M\times \R^{|\mathcal{S}|}, \lbb I, \ed u^\alpha - \lambda^\alpha \ed t\rbb_{\alpha\in \mathcal{S}}; \ed t),
		\]
	where $\mathcal{S} \subseteq \{1,2,\ldots,m\}$ is fixed, is called a \emph{prolongation by differentiation} of $(M,I;\tau)$; note that this is a particular type of partial prolongation.}
	
	\item{A \emph{Cartan prolongation} of $(M,I;\tau)$ is a system $(N,J;\sigma)$ together with a submersion
		$\pi: N\rightarrow M$ satisfying
		\begin{enumerate}[\bf i.]
			\item{$\pi^*I\subset J;$}
			\item{$\pi^*\tau = \sigma$;}
			\item{any (generic) integral curve\footnote{\label{CPFN}
			By ``generic," we mean that this condition holds for a set of integral curves that is open and dense in the space of all integral curves with respect to the $C^k$ topology for some $k$.
We do \emph{not} require all integral curves of 
				$(M,I;\tau)$ to satisfy the condition, in order to make the concept useful even when
				there exists a negligible set of integral curves that do not admit unique liftings. For example,
				consider $(M,I;\tau) := (\R^5, \lbb\ed x - u \ed t, \ed y - v\ed t\rbb; \ed t)$.
	  			It has a Cartan prolongation: 
				$(N,J;\sigma) := (\R^5\times\R, \lbb I, \ed u - \lambda \ed v\rbb; \ed t)$, to which all
				integral curves of $(M,I;\tau)$ admits a unique lifting except those along which 
				$v$ is a constant.} $\gamma: (-\epsilon,\epsilon)\rightarrow M$ of $(M,I;\tau)$
				has a unique lifting $\hat\gamma:(-\epsilon,\epsilon)\rightarrow N$ that satisfies
					$\hat\gamma^*J = 0$ and $\pi\circ\hat\gamma = \gamma$;}
		\end{enumerate}
		The fiber dimension of $\pi$ will be called the \emph{order} of the Cartan prolongation.
		}
	\end{enumerate}
	\end{Def}	

%%Comment
%{\color{red} 3. These next two paragraphs were pretty confusing to me on a first read-through; can you have a go at making them more clear?\bf (Response: These paragraphs are revised. Proposition 2.1 added.} {\color{blue} Comment: This is very helpful!  I made some minor edits.)}

\begin{remark}\label{cauchyCharRemark}
	One can ask:
	\emph{Does a partial prolongation of a system necessarily yield a system?} 
	This is the question that 
	motivated us to include the condition {\bf ii} in Definition \ref{systemDef}. Without this condition, 
	if a Pfaffian system $(M,I)$ with an independence condition $\tau$ 
	admits an integral surface $\iota: S\hookrightarrow M$ such that $\iota^*\tau$ is nonvanishing, 
	then a partial prolongation may 
	introduce nontrivial Cauchy characteristics, which is undesirable. As an example, consider the system
		\[
			I  = \lbb\ed z - p \ed x - q\ed t\rbb,
		\]
	defined on $M = \R^5$ with coordinates $(x,t,z,p,q)$ and with independence condition $\ed t$. The Pfaffian system $(M,I)$ admits integral surfaces.
	A partial prolongation may be achieved by 
	adjoining to $I$
	the $1$-form
		\[
			\ed x - \lambda \ed t.
		\]
	The resulting system has nontrivial Cauchy characteristics, corresponding to those curves that annihilate the $1$-forms below:
		\[
			\ed z, ~ \ed (\lambda p+ q), ~ \ed x, ~\ed \lambda, ~\ed t.
		\]	
	On the other hand, we have the following Proposition.	
	\end{remark}

\begin{Prop}
	A partial prolongation of a system remains a system.
\end{Prop}	
	
\emph{Proof.} 
Let $(M,I;\tau)$ be a system. Let $(N,J;\tau)$ be a partial prolongation of $(M,I;\tau)$, obtained by adjoining 
the $1$-forms  (as in Definition \ref{ProlongationsDef})
	\[
		\mu^\alpha - \lambda^\alpha \tau, \qquad \alpha \in \mathcal{S}\subseteq \{1,2,\ldots,m\}.
	\]
To show that $(N,J;\tau)$ is a system, it suffices to verify that
 (1) it does not admit any Cauchy characteristics and (2)
it does not admit any integral surfaces on which $\tau$ is nonvanishing.	

If $(N,J;\tau)$ has Cauchy characteristics $\mathcal{C}$, then any integral curve transversal to
a curve in $\mathcal{C}$ and satisfying the independence condition $\tau$
can be used to generate an integral surface on which $\tau$ pulls back to be nonzero.
Thus, it suffices to justify (2). 

Suppose that $\iota:S\rightarrow N$ is an integral surface of $(N,J)$ satisfying $\iota^*\tau \ne 0$.
Let $\pi: N\rightarrow M$ be the obvious submersion. Since $\pi^*I\subset J$ and $\iota^*\tau \ne 0$,
locally the rank of $\pi|_S$ must be equal to 1; otherwise, $(M,I)$ would admit an integral surface where $\tau$
pulls back to be nonzero. As a result, $\iota^*\mu^\alpha$ are multiples of $\tau$, and $\iota^*\ed \mu^\alpha$ are equal to zero.
Since 
	\[
		0=\iota^*\ed(\mu^\alpha - \lambda^\alpha\tau) = \iota^*(\ed \mu^\alpha - \ed \lambda^\alpha\W\tau), \quad \alpha\in \mathcal{S},
	\]
$\ed \lambda^\alpha$ must all be multiples of $\tau$; in other words, if locally $\tau = \ed t$, then 
$\lambda^\alpha$ are functions of $t$,  which violates the assumption that $S$ is a surface.\qed

{\remark
	We note that each notion of prolongation in Definition \ref{ProlongationsDef}
	requires that the independence conditions
	correspond to each other via the
	underlying submersion.  Removing this requirement leads to broader notions of prolongation
	that are more familiar in the literature. In this paper, we will always assume that independence 
	conditions are matched by a prolongation.
	
}	
	
{\remark	\label{prolongationRmk} The notions of prolongation above satisfy the following.\footnote{The reader
may compare this with
		\cite[p.61]{Sluis}.}
	\begin{enumerate}[\bf A.]
		\item{Among the four notions above, Cartan prolongation is the most general.}
		\item{An order-$1$ Cartan prolongation
			of a system $(M,I;\tau)$ need not be a partial prolongation of $(M,I;\tau)$.
			(See \cite[p.50]{Sluis}.)
			}

		\item{If $(M,I;\tau)$ is a CTS, then a prolongation by differentiation of $(M,I;\tau)$ is also 
			a CTS.}
		\item{A partial prolongation may not necessarily be realized as a prolongation by differentiation of the 
		same system. For example, let $M = \R^7$ and $I = \lbb \theta^1,\theta^2,\theta^3\rbb$, where
		$\theta^i = \ed x^i - p^i\ed t$ $(i = 1,2,3)$. The system
			\[
				(N,J;\tau):= (\R^7\times\R, \lbb I, \ed p^1 + p^2 \ed p^3 - \lambda \ed t\rbb; \ed t)
			\]
		is a partial prolongation of $(M,I;\ed t)$. However, it is not a prolongation by differentiation,
		since $\lbb J, \ed t\rbb$ is not Frobenius (in other words, $(N,J;\tau)$ is not a CTS).
		\item{When $m = \dim(M) -\rank(I) - 1 = 1$, a Cartan prolongation is necessarily the
		result of successive total prolongations. This case is treated in \cite{Cartan14}, where Cartan 
		did not require the matching of independence conditions.}
	}
	    \item{A composition of Cartan prolongations is also a Cartan prolongation, as the following proposition shows.}
	
	\end{enumerate}
}

		\begin{Prop}
			If $\pi: (N,J;\sigma)\rightarrow (M,I;\tau)$ and $\varpi: (P,L; \rho)\rightarrow(N,J;\sigma)$
			 are both Cartan prolongations of systems, then the composition
			$\pi\circ\varpi: (P,L;\rho)\rightarrow(M,I;\tau)$ is again a Cartan prolongation.
		\end{Prop}
		\emph{Proof.} By the assumption,
					\[
						(\pi\circ\varpi)^* I  = \varpi^*(\pi^*I)\subseteq \varpi^*J\subseteq L,
					\]
				and it is clear that each generic integral curve $\gamma$ of $(M,I;\tau)$ has a lifting to
				$(N,J;\sigma)$ then to $(P,L;\rho)$. To justify uniqueness, suppose that 
				$\gamma_1$ and $\gamma_2$ are two liftings of $\gamma$ to $(P,L;\rho)$.
				Therefore both $\varpi\circ\gamma_i$ $(i = 1,2)$ are
				integral curves of $(N,J;\sigma)$ and project via $\pi$ to $\gamma$. 
				Since $\pi$ is a Cartan prolongation,
				$\varpi\circ\gamma_1 = \varpi\circ\gamma_2$; since $\varpi$ is a Cartan prolongation,
				$\gamma_1 = \gamma_2$.\qed

	\subsection{The Extension Theorem of Sluis}\label{extThmsec}
	
		In \cite{Cartan14}, Cartan noted that a generalization of the fact {\bf E} in 
		Remark \ref{prolongationRmk}
		into the cases of $m\ge 2$ can be quite difficult. 
		
		Sluis, in his thesis \cite{Sluis}, considered 
		such more general cases. In this section, we remind the reader of a notable theorem 
		he obtained and sketch
		the main arguments in his proof. The theorem indicates the following. (See Theorem~\ref{SluisThm} for a more
		precise statement.)
		\begin{quote}
		\emph{Any Cartan prolongation of a system can be extended to a (successive) total prolongation, and such an
		extension itself is also a Cartan prolongation. }
		\end{quote}
		
		Consider a Cartan prolongation of order $r$ represented by the diagram below,
		where the rank of $I$ is $n$.

		{\begin{figure*}[h!]
					\begin{center}
					\begin{tikzcd}[column sep= 0.1mm]
		(N^{n+m+r+1}, J;\sigma)\arrow[d,swap,"\displaystyle\pi"]  \\[0.1em]
						(M^{n+m+1},I;\tau)                 										
 					\end{tikzcd}
					\end{center}
				\end{figure*}}
				
		By shrinking $M$ if needed, choose coordinates $(x^i,u^\alpha,t)$ on $M$ such that 
		$I, \ed u^\alpha, \ed t$ span $T^*M$, where $\ed t = \tau$.
		
%%Comment		
%{\color{red}  4. Say something about WHY it might occur to one to ask this!  Also, note that you haven't defined $I_k$ yet - this doesn't come until Definition 2.8 - so introducing it here with a bracket around it seems poorly motivated.  Would it make sense to move that definition earlier, or at least mention here that $I_1$ is the first of a sequence of bundles that will be defined later? \bf(Response: The confusing notation is removed; bundle notation replaced by $\hat J$. See the paragraph below.)} {\color{blue} Oh, this is much better - thanks!}

		In order to understand such a Cartan prolongation, Sluis began by considering
		those $1$-forms in $J$ that are
		closest to being expressible in terms of the coordinates on $M$. To be more specific, he asked:
		\emph{What is the rank of the subbundle $\hat J$ of $J$
		spanned by the $1$-forms that can be written as
				\[
					f_1\ed u^1+f_2\ed u^2+\cdots + f_m \ed u^m + g \ed t,
				\]
		where $f_\alpha, g$ are functions on $N$?}	(Note that $\hat J$ has the coordinate-free interpretation as a quotient bundle $(J\cap \pi^*(T^*M))/ \pi^*I$.)

We have the following lemma:				
\begin{lemma}\label{dim-J-hat-lemma}		
			$\rank(\hat J)\ne 0, m+ 1$.
\end{lemma}

\begin{proof}
			Let $(x^i, u^\alpha, v^\rho, t)$ $(\rho = 1,\ldots,r)$ 
			be local coordinates on $N$.
			
			If $\rank(\hat  J) = 0$, then 
			$J$ is spanned by $I$ and some $1$-forms
				\[
					a^k_\rho \ed v^\rho  + b^k_\alpha \ed u^\alpha + c^k\ed t 
					\quad (k = 1,\ldots, s;~ s\le r),
				\]
			where $a^k_\rho,b^k_\alpha,c^k$ are functions on $N$, and
			$(a^k_\rho)$ has maximum rank (i.e., rank $s$). 
			Now let $\gamma = (\x(t), \u(t), t)$ be any integral curve of 
			$(M,I;\tau)$. The pullback of $J$ to $\pi^{-1}\gamma$ is spanned by the $1$-forms
				\[
					a^k_\rho \ed v^\rho + (b^k_\alpha (u^\alpha)' + c^k)\ed t.
				\]
			Since the row rank of $(a^k_\rho)$ is full,  $J$ induces a distribution $J^\perp$ on $\pi^{-1}\gamma$, 
			whose integral curves are non-unique  
			liftings of $\gamma$. This violates part {\bf iii} in the definition of Cartan prolongation.
				
			If $\rank(\hat J)=m+1$, then $J$ would contain the independence condition,
			which violates the definition of a system.
\end{proof}

		By Lemma \ref{dim-J-hat-lemma}, the only cases that can occur are:
		\[
		\mbox{\textbf{I.} } \rank(\hat J) = m\qquad and \qquad 
		\mbox{\textbf{II.} }  1\le \rank(\hat J)<m.
		\]
			
		\noindent{\bf Case I.} 
		In this case, there exist functions $ f^\alpha: N\rightarrow \R$ $(\alpha = 1,\ldots,m)$
		such that 
		\begin{equation}\label{CartanProCaseI_fdefn}
			\ed  u^\alpha  -  f^\alpha \ed t
		\end{equation}
		are sections of $J$.
		It is shown in \cite[p.65]{Sluis} that, on
	some dense open subset of $N$, $f^\alpha$ must be  	
		independent of $x^i, u^\alpha, t$ and independent among themselves; 
		otherwise, $(f^\alpha)$ cannot be surjective for a 
	fixed initial point on $M$ (by Sard's theorem),
	and a generic integral curve $\gamma$ of $(M,I;\tau)$ passing through that
	 initial point would not admit a lifting.
	
		Thus, by shrinking $N$, if needed, one can define a submersion:
		\[
			\pi^1: N\rightarrow M\times \R^m
		\]
		by 
		\[
			N\ni p\mapsto (\pi(p), (f^\alpha(p))).
		\]
	We can identify $M\times \R^m$ with the space of the total prolongation $\pr^{(1)}M$
	equipped with the Pfaffian system
	generated by $\ed u^\alpha - y^\alpha \ed t$ $(\alpha = 1,\ldots,m)$, where $(y^\alpha)$ are coordinates
	on the $\R^m$-component. This makes $\pi^1$ a Cartan prolongation.	
	See the diagram below.\footnote{Here and below, 
	we drop the independence conditions in these
	diagrams for clarity.}
	{\begin{figure*}[h!]
					\begin{center}
					\begin{tikzcd}[column sep = 5, row sep = 8]
		(N, J)\arrow[dd,swap,"\displaystyle\pi"] \arrow[rd, "\displaystyle\pi^1"] \\[0.1em]
				& (\pr^{(1)}M,\pr^{(1)}I)\arrow[ld,"\displaystyle\pi_1"]\\[0.1em]
						(M,I)                 										
 					\end{tikzcd}
					\end{center}
				\end{figure*}	
	} 

	Now note that $(\pr^{(1)}M, \pr^{(1)}I;\ed t)$ has type $(n+m, m)$.
	$\pi^1$ is a Cartan prolongation of order $(r - m)$.\\

	\noindent{\bf Case II.} In this case, $\rank(\hat J) = q<m$.
	Suppose that $\hat J$ has the following basis representatives:
	\begin{align*}
		f^\mu_1\ed u^1 + &\cdots+ f^\mu_m \ed u^m + g^\mu \ed t\quad (\mu = 1,\ldots,q).
	\end{align*}
	Let $\mu,\nu = 1,\ldots,q$ and $\alpha = 1,\ldots,m$.
	We must have $\rank(f^{\nu}_\alpha) = q$, since, otherwise, $\ed t$ would be a section of
	 $J$.	

	By reordering $\ed u^\alpha$, we may assume that
		\[
			\det(f^\nu_\mu)  \ne 0.
		\]
	Now consider the following prolongation of $(N,J)$ by differentiation:
	\[
		N_1:= N\times \R^{m-q},
	\]
	with $(y^1,\ldots,y^{m-q})$ being coordinates on the $\R^{m-q}$-component, and
	\[
		J_1 := \lbb J, \ed u^{q+1} - y^1\ed t,\ldots, \ed u^{m} - y^{m-q}\ed t\rbb.
	\]
	Let the submersion $N_1\rightarrow N$ be denoted by $\phi$.
	
	It is easy to see that \[\pi\circ\phi: (N_1, J_1)\rightarrow (M,I)\] represents
	 a Cartan prolongation that
	belongs to {Case \bf I}.
	This implies the diagram in Figure~\ref{Fig_CartanCaseII}, where $\pi^1$ is a Cartan prolongation of order $r - q$:
	{\begin{figure}[h!]
					\begin{center}
					\begin{tikzcd}[column sep= 18, row sep = 15]
		(N, J)\arrow[d,swap,"\displaystyle\pi"]& (N_1, J_1) \arrow[d, "\displaystyle\pi^1"]\arrow[l,swap, "\displaystyle \phi"]\arrow[ld] \\[0.2em]
		(M,I)		& (\pr^{(1)}M,\pr^{(1)}I)\arrow[l,"\displaystyle\pi_1"]
 					\end{tikzcd}
					\end{center}
		\caption{Extending a Cartan prolongation in Case \textbf{II}.}
		\label{Fig_CartanCaseII}
	\end{figure}	
	} 
	
	Combining the cases {\bf I} and {\bf II}, for some minimal finite $K$, we obtain the diagram in Figure~\ref{Fig_SluisExt}, 
		where $\pi^K$ is an isomorphism.
{\begin{figure}[h!]
					\begin{center}
					\begin{tikzcd}[column sep= 18, row sep = 15]
	(N, J)\arrow[d,swap,"\displaystyle\pi"]& (N_1, J_1) \arrow[d, "\displaystyle\pi^1"]\arrow[l, swap,"\displaystyle \phi_1"] &\cdots\arrow[l, swap,"\displaystyle \phi_2"] 
	& (N_K, J_K) \arrow[d,equal ,"\displaystyle\pi^K"]\arrow[l, swap,"\displaystyle \phi_K"]\\[0.2em]
	(M,I)		& (\pr^{(1)}M,\pr^{(1)}I)\arrow[l,"\displaystyle\pi_1"]&
	 \cdots\arrow[l,"\displaystyle\pi_2"]
	 & (\pr^{(K)}M,\pr^{(K)}I)\arrow[l,"\displaystyle\pi_K"]
 					\end{tikzcd}
					\end{center}
					\caption{Successive extensions of Cartan prolongations.}
					\label{Fig_SluisExt}
	\end{figure}	
	} 
	
	In Figure~\ref{Fig_SluisExt}, 
	\begin{enumerate}[(1)]
		\item{ Each $\pi_k$ has order $m$;}
		\item{Each $\phi_k$ has order $m - q_k$, where $1\le q_k\le m$. When $q_k = m$, $\phi_k$
				is constructed from Case {\bf I} and is an isomorphism;}
		\item{Each $\pi^k$ has order $r - q_1 - q_2-\cdots - q_k$.}
		\item{$r = q_1+\cdots + q_K$, and consequently 
					$\left\lfloor \frac{r-1}{m}\right\rfloor + 1\le K \le r$.}
	\end{enumerate}	

	The construction above is summarized by the following {extension theorem} of Sluis; it will be illustrated in Examples \ref{example-4} and \ref{example-5} below.
	\begin{theorem}\cite[p.63, Theorem 24]{Sluis}
	{If $\pi:(N,J;\sigma)\rightarrow (M,I;\tau)$ is a Cartan prolongation of order $r$, 
	then there exists an integer $K\le r$ and a map $\hat \pi$ given by a composition of $K$ successive prolongations by differentiation such that
		the following diagram commutes in the sense of EDS, where $\pi_{K,0}: \pr^{(K)}M\rightarrow M$ is the composition of $K$ successive total prolongations.}
		{\begin{figure*}[h!]
					\begin{center}
					\begin{tikzcd}[column sep = 5, row sep = 8]
		(\pr^{(K)}M,\pr^{(K)}I)\arrow[dd,swap,"\displaystyle\pi_{K,0}"] \arrow[rd,"\displaystyle\hat\pi"]	& \\[0.1em]
			 	     &  (N,J)\arrow[ld,"\displaystyle\pi"]	\\
						 		(M,I) &                										
 					\end{tikzcd}
					\end{center}
				\end{figure*}	}	\label{SluisThm}
	\end{theorem}

	 \begin{corollary}\label{simpleCPrankCond}
		 	If $\pi:(N,J;\sigma)\rightarrow(M,I;\tau)$ is a Cartan prolongation, then
			\[
				\rank(J) - \rank(I) =  \dim N - \dim M.
			\]
		 \end{corollary}
		
	\emph{Proof.} This is because, in the argument that leads to Theorem \ref{SluisThm},
			\begin{equation*}
			\begin{alignedat}{2}
				\dim N_k - \dim N_{k-1} &= \rank (J_k)   - \rank (J_{k - 1}),\\
				\dim \pr^{(k)}M - \dim \pr^{(k-1)}M &= \rank (\pr^{(k)}I) - \rank  (\pr^{(k-1)}I),
			\end{alignedat}
			\end{equation*}
			and  at the $K$-th stage, $(N_K,J_K)$ and $(\pr^{(K)}M,\pr^{(K)}I)$ coincide.
			\qed

	{\remark
		At various points in the argument above, we have applied steps such as 
		``by shrinking to an open dense subdomain, if needed \dots". This is because, as Footnote \ref{CPFN}
		indicated, we allow a negligible set of integral curves to ill-behave relating to
		a Cartan prolongation. Putting the example in Footnote \ref{CPFN} in context, one 
		immediately notices that that Cartan prolongation belongs to Case {\bf II}, and we have
		\begin{align*}
			(N_1, J_1)&= (N\times \R, \lbb J, \ed v - \mu \ed t\rbb) \\
					&= (N\times \R, \lbb \ed x - u\ed t, 
						\ed y - v\ed t, \ed u - \lambda \ed v, \ed v - \mu \ed t\rbb),
		\end{align*}
		where $\mu$ is the coordinate on the $\R$-component of $N_1$. 
		Thus,
		 the system $(N_1, J_1)$, 
		relative to $(M, I)$, is 
		a Cartan prolongation in Case {\bf I}. In particular, 
			\[
				\ed u - \lambda\mu\ed t, \quad \ed v -\mu \ed t
			\]
		are sections of $J_1$. According to \eqref{CartanProCaseI_fdefn},
			\[
				f^1 = \lambda\mu, \quad f^2 = \mu.
			\]
		Note that $\ed f^1, \ed f^2$ are only linearly independent when $\mu\ne 0$,
		so
			\[
				\pi^1: (N_1, J_1) \xrightarrow{\cong} (\pr^{(1)}M, \pr^{(1)}I), 
			\]
		 is understood as a diffeomorphism onto a dense open subset of $\pr^{(1)}M$. 
			}

	\subsection{Absolute $\Leftrightarrow$ Dynamic}

	A consequence of Sluis's Extension Theorem is that \emph{absolute equivalence} 
	is an equivalence relation; this
	is proved in \cite{Sluis}. 
	Moreover, under some mild assumptions,  absolute equivalence is equivalent 
	to the notion of \emph{dynamic equivalence}, as we will
	demonstrate below.
		
	\begin{Def}\label{AbsEqDef}
		Two systems $(M,I;\tau)$ and $(\bar M, \bar I;\bar\tau)$ are said to be \emph{$\ut$-absolutely equivalent} if there exists a system $(N,J;\sigma)$
	and submersions $\pi: N\rightarrow M$ and $\bar \pi:N\rightarrow \bar M$ that realize $(N,J;\sigma)$
	as a Cartan prolongation of both $(M,I;\tau)$ and $(\bar M,\bar I;\bar\tau)$.
	{\begin{figure*}[h!]
					\begin{center}
					\begin{tikzcd}[column sep= 0.2mm]
	&	(N,J;\sigma)\arrow[dl,swap,"\displaystyle\pi"] \arrow[rd,"\displaystyle\bar\pi"]	& \\
			(M,I;\tau) 	     &  &     (\bar M,\bar I;\bar\tau)
 					\end{tikzcd}
					\end{center}
		\end{figure*}	
	}
	\end{Def}	
		
	\begin{Def}\label{DynEqDef}
		Two systems $(M,I;\tau)$ and $(\bar M,\bar I;\bar\tau)$ are said to be 
		\emph{$\ut$-dynamically equivalent}\footnote{Compare with \cite[Definition 3.2]{NRM94}.} if
	there exist integers $p,q\ge 0$ and submersions $\Phi,\Psi$, as shown in the diagram below, such that
	\begin{enumerate}[\bf i.]
	\item{$\Phi^*\bar\tau  = \pi^*\tau$, $\Psi^*\tau = \bar\pi^*\bar\tau$;}
	\item{$
		\Phi^* \bar I\subset \pr^{(p)}I, ~\Psi^*  I\subset \pr^{(q)}\bar I
	$}
	\item{for any (generic) integral curve $\gamma$ of $(M,I;\tau)$, 
	$\Psi\circ(\Phi\circ\gamma^{(p)})^{(q)} = \gamma$,
	and for any (generic) integral curve $\bar\gamma$ of $(\bar M,\bar I;\bar\tau)$, $\Phi\circ(\Psi\circ\bar\gamma^{(q)})^{(p)} = \bar\gamma$.}
	\end{enumerate}
\begin{figure*}[h!]
	\begin{center}
		\begin{tikzcd}[column sep = 25, row sep = 25]
		&   \pr^{(p)} M \arrow[dr,"\phantom{aa}\Phi", sloped]\arrow[d,"\pi", swap]  & 	 \pr^{(q)}\bar M \arrow[dl, "\Psi\phantom{aa}",swap, sloped]\arrow [d,"\bar\pi"]	&\\
		\R \arrow[ur,"\gamma^{(p)}"]\arrow[r,"\gamma"]& M & \bar M & \R \arrow[ul,"\bar\gamma^{(q)}", swap]\arrow[l,"\bar\gamma", swap]
 					\end{tikzcd}
				\end{center}	
\end{figure*}
	\end{Def}	
	
	\begin{theorem}\label{AbsDyn}
		Two systems $(M,I;\tau)$ and $(\bar M, \bar I;\bar \tau)$ are $\ut$-absolutely equivalent
		if and only if they are $\ut$-dynamically equivalent.\footnote{This theorem may be seen as a variant of \cite[Theorem 3.6]{NRM94}}
	\end{theorem}	
	\emph{Proof.}
	 $(\Rightarrow)$ Start with a $\ut$-absolute equivalence as described in Definition \ref{AbsEqDef}.
	 	By Theorem \ref{SluisThm}, there exist integers 
		$p,q$ such that the following diagram commutes, providing a bijection between 
		(generic) integral curves of each system involved. It follows by definition 
	that $(M,I;\tau)$ and $(\bar M, \bar I;\bar \tau)$ are $\ut$-dynamically equivalent.
		{\begin{figure*}[h!]
			\begin{center}
					\begin{tikzcd}[column sep= 5, row sep = 8]
					
	(\pr^{(p)}M,\pr^{(p)}I)\arrow[dr,swap,"\displaystyle \phi"] \arrow[dd,swap,"\displaystyle\pi_{p,0}"]&& 	(\pr^{(q)}\bar M,\pr^{(q)}\bar I)\arrow[dl,"\displaystyle \psi"]\arrow[dd,"\displaystyle\bar\pi_{q,0}"]\\			
	&	(N,J)\arrow[dl,"\displaystyle\pi"] \arrow[rd,swap,"\displaystyle\bar\pi"]	& \\
			(M,I) 	     &  &     (\bar M,\bar I)     										
 					\end{tikzcd}
					\end{center}
				\end{figure*}	
	}

	 $(\Leftarrow)$ Conversely, assume the diagram in Definition  \ref{DynEqDef}. It suffices to show
	 that $\Psi$ represents a Cartan prolongation. 

Let $\gamma: (-\epsilon,\epsilon)\rightarrow M$ be a generic integral curve of $I$.  By the third condition in Definition \ref{DynEqDef}, the curve 
\[ (\Phi\circ\gamma^{(p)})^{(q)}: (-\epsilon,\epsilon)\rightarrow \pr^{(q)} \bar M \]
is a lifting of $\gamma$ into $\pr^{(q)} \bar M$ as an integral curve of $\pr^{(q)} I$.

	 Now suppose that there are two such liftings 
	 \[\gamma_1,\gamma_2: (-\epsilon,\epsilon)\rightarrow \pr^{(q)}\bar M. \] 
 By construction, $\bar\gamma_1:= \bar\pi\circ\gamma_1$ and $\bar\gamma_2:= \bar\pi\circ\gamma_2$ satisfy
	 	\[
			\Psi\circ\bar\gamma_1^{(q)} = \Psi\circ\bar\gamma_2^{(q)}.
		\]
	Hence,
		\[
			\bar\gamma_1 = \Phi\circ\(\Psi\circ \bar\gamma_1^{(q)}\)^{(p)}
				 = \Phi\circ\(\Psi\circ\bar\gamma_2^{(q)}\)^{(p)} = \bar\gamma_2.
		\]		
	 Since $\bar\pi$ is a total prolongation, it follows that 
	 	\[
			\gamma_1 = \gamma_2.
		\]
	This proves that $\Psi$ is a Cartan prolongation. 
	
	We end the proof by remarking that 
	the independence conditions are preserved by all the maps involved.
	\qed

	\subsection{Relative Extensions}\label{relativeExtsec}
		
	In order to understand Cartan prolongations further, we introduce the notion of 
	a \emph{relative extension}.
	
	\begin{Def}\label{relExtDef}
		Let $\pi: (N,J;\sigma)\rightarrow (M,I;\tau)$ be a Cartan prolongation. 
		We define the \emph{$k$-th extension $I_k$ of $I$ relative to $\pi$} inductively:
			\begin{enumerate}[(1)]
				\item{$I_0 = \pi^*I$;}
				\item{$I_k = \mathcal{C}(I_{k-1})\cap J$\quad  $(k\ge 1)$.}
			\end{enumerate}
	\end{Def}

%%Comment	
%{\color{red} 6. I think a simple example here would be instructive - it could be VERY simple, even something like a Brunovsk\`y normal form. \bf (Response: Example added below.)}	{\color{blue} Great - this is exactly what I had in mind!  I did some minor reformatting and editing.} {\color{orange} \bf Thank you!}
	
	\begin{Example}\label{relExtExample}
		Let $\theta^i_j$ denote the $1$-forms
			\[
				\theta^i_j = \ed x^i_j - x^i_{j+1}\ed t.
			\]
		Let $(M,I;\ed t)$ be the type $(3,2)$ CTS in Brunovsk\'y normal form generated by the three $1$-forms
\begin{alignat*}{2}
\theta^1_0 & = \ed x^1_0 - x^1_1 \ed t,& \qquad \theta^2_0 &= \ed x^2_0 - x^2_1 \ed t, \\
\theta^1_1 & = \ed x^1_1 - x^1_2 \ed t. & &
\end{alignat*}

		This system has a prolongation by differentiation to the type (6,2) system $(N,J;\ed t)$ in Brunovsk\'y normal form generated by the six $1$-forms
\begin{alignat*}{2}
\theta^1_0 & = \ed x^1_0 - x^1_1 \ed t, & \qquad \theta^2_0 &= \ed x^2_0 - x^2_1 \ed t, \\
\theta^1_1 & = \ed x^1_1 - x^1_2 \ed t. & \qquad \theta^2_1 &= \ed x^2_1 - x^2_2 \ed t, \\
\theta^1_2 & = \ed x^1_2 - x^1_3 \ed t, & & \\
\theta^1_3 & = \ed x^1_3 - x^1_4 \ed t. & &
\end{alignat*}
		The relative extensions of $I$ are
			\[
				I_0 = \pi^*I,\quad I_1 = \lbb I, \theta^1_2, \theta^2_1\rbb, \quad I_2 = J.
			\]	
	\end{Example}

		Definition \ref{relExtDef} is independent of the choice of coordinates. 
		Moreover, 
		there exists an integer $K\ge 0$ 
		indicating where $I_k$ stabilizes: 
		\[
			\pi^*I = I_0\subsetneq I_1\subsetneq\cdots\subsetneq I_K = I_{K+1} = \cdots \subseteq J.
		\] 	
		For simplicity,  we denote $I_\infty:= I_K$, and we define the \emph{extension length}
		of $\pi$ to be the smallest integer $K$ satisfying $I_\infty = I_K$.
		
		 Each $I_k\subseteq J$ is a subbundle on $N$,
			so they may admit nontrivial Cauchy characteristics.
			As a result, we consider the underlying manifold of $I_k$ as the one
			determined by
		 the Cartan system $\mathcal{C}(I_k)$. 
			However, since an inclusion between Pfaffian bundles does not
			generally imply the inclusion of their Cartan systems (for instance, 
			$\lbb\ed y - z\ed x\rbb\subsetneq\lbb\ed x, \ed y\rbb$, but their Cartan systems
			satisfy $\lbb \ed x,\ed y,\ed z\rbb\supsetneq\lbb\ed x,\ed y\rbb$), it is not immediately
			clear how the $\mathcal{C}(I_k)$ relate to each other 
			or whether $I_\infty$ must be equal 
			to $J$ in general.

	\subsection{Regular and $\mathcal{C}$-regular Cartan Prolongations}

		\begin{Def}
			We call a Cartan prolongation $\pi:(N,J;\sigma)\rightarrow(M,I;\tau)$ \emph{simple}
			if $I_1 = J$.
		\end{Def}
		
		{\remark Because a system is assumed to have no Cauchy characteristics, it is easy to see that the total prolongation, any partial prolongation or a prolongation by differentiation
		 (Section \ref{ProlongationSection}) of a system
		are simple Cartan prolongations.}
		
		\begin{lemma}\label{rank1CPsimple}
			Any Cartan prolongation of order $1$ is simple.
		\end{lemma}
		
		\emph{Proof.} In this case, if $I_1\ne J$, then $I_1 = \pi^*I$, which 
				is impossible by Lemma \ref{dim-J-hat-lemma}. \qed\\

		The following theorem, which will be useful later,
		can be regarded as a special case of Theorem \ref{SluisThm}.
		
		\begin{theorem}\label{simpleCPThm}
			If $\pi:(N,J;\sigma)\rightarrow(M,I;\tau)$ is a simple Cartan prolongation, then there exists
			a simple Cartan prolongation 
				$\hat\pi: (\pr^{(1)}M,\pr^{(1)}I)\rightarrow (N,J)$ such that the diagram 
				below commutes in the EDS sense.
					{\begin{figure*}[h!]
					\begin{center}
					\begin{tikzcd}[column sep = 5, row sep = 8]
		(\pr^{(1)}M,\pr^{(1)}I)\arrow[dd,swap,"\displaystyle\pi_{1,0}"] \arrow[rd,"\displaystyle\hat\pi"]	& \\[0.1em]
			 	     &  (N,J)\arrow[ld,"\displaystyle\pi"]	\\
						 		(M,I) &                										
 					\end{tikzcd}
					\end{center}
				\end{figure*}	}
		\end{theorem}
		
		\begin{remark}\label{simpleCPrmk}
			$\hat\pi$ is in fact a prolongation by differentiation.
		\end{remark}
			 
		With Theorem \ref{simpleCPThm} in mind, we are interested in the case when a Cartan prolongation
		can be achieved by successively performing simple Cartan prolongations, starting from an original
		system. To make this point explicit, we make the definition below.

		\begin{Def}\label{regularity}
			Let $\pi: (N,J;\sigma)\rightarrow(M,I;\tau)$ be a Cartan prolongation.
			It is
			\emph{regular} if there exist vector subbundles $I_{(\ell)}\subset J$ $(\ell  = 1,\ldots, L)$ satisfying
				\[
					\pi^*I  = I_{(0)}\subsetneq I_{(1)}\subsetneq\cdots \subsetneq I_{(L-1)}\subsetneq I_{(L)} = J
				\]
			such that
			\begin{enumerate}[(1)]
				\item{$\sigma$ is a section of $\mathcal{C}(I_{(\ell)})$ for each $\ell\in \{0,1,\ldots, L\}$;}
				\item{each $(M_{(\ell)}, I_{(\ell)};\sigma_{(\ell)})$ is a system, where $M_{(\ell)}$ is the manifold given by the quotient space of $N$ by the leaves of the distribution annihilated by $\mathcal{C}(I_{(\ell)})$, and $\sigma_{(\ell)}$ is a corresponding independence condition defined
				 on $M_{(\ell)}$;} 
				\item{each inclusion $I_{(\ell)}\subsetneq I_{(\ell+1)}$ induces an inclusion $\mathcal{C}(I_{(\ell)}) \subsetneq \mathcal{C}(I_{(\ell+1)})$, which in turn determines a submersion
				from $M_{(\ell+1)}$ to $M_{(\ell)}$ that
				represents a simple Cartan prolongation 
				of $(M_{(\ell)}, I_{(\ell)};\sigma_{(\ell)})$.}
			\end{enumerate}	
			Otherwise, $\pi$ is called \emph{singular}.
		\end{Def}

		A condition that is stronger than ``regular" is when the vector subbundles in 
		Definition \ref{regularity} can be chosen to be the canonical relative extensions $I_k$ 
		(Definition \ref{relExtDef}) and still satisfy the conditions (1)-(3). To be clear, we present the
		 following definition.

		\begin{Def}\label{Cregularity}
			Let $\pi: (N,J;\sigma)\rightarrow(M,I;\tau)$ be a Cartan prolongation.
			It is $\mathcal{C}$-\emph{regular} if the canonical relative extensions $I_k$ ($k = 1,\ldots K$) (see Definition \ref{relExtDef}) satisfy
			\[
				\pi^*I = I_0 \subsetneq I_1\subsetneq\cdots\subsetneq I_{K-1}\subsetneq I_K = J,
			\]
			and 
			\begin{enumerate}[(1)]
				\item{$\sigma$ is a section of $\mathcal{C}(I_k)$ for each 
							$k\in \{0,1,\ldots, K\}$;}
				\item{each $(M_k, I_k;\sigma_k)$ is a system, where
					$M_k$ stands for the manifold determined by $\mathcal{C}(I_k)$,
					and $\sigma_k$ is a corresponding independence condition defined
					on $M_k$;}
				\item{each inclusion $I_k\subsetneq I_{k+1}$ induces a submersion 
				from $M_{k+1}$ to $M_{k}$, which represents a Cartan prolongation
				of $(M_k, I_k;\sigma_k)$.}
			\end{enumerate}
		\end{Def}
		
		\begin{Prop}
			If $\pi: (N,J;\sigma)\rightarrow (M,I;\tau)$ is $\mathcal{C}$-regular, then 
			the extensions $I_k$ of $I$ relative to $\pi$ satisfy the condition that $I_k\subsetneq I_{k+1}$ is a 
			simple Cartan prolongation.
		\end{Prop}
		\emph{Proof.} It suffices to prove that the Cartan prolongation induced by
		 $I_k\subsetneq I_{k+1}$ is simple.
				This is immediate since
					$
						I_{k+1} = \mathcal{C}(I_k)\cap J = \mathcal{C}(I_k)\cap I_{k+1}.
					$
					\qed
		
		\begin{lemma}\label{regFilVSrelExt}
			Let $\pi:(N,J;\sigma)\rightarrow (M,I;\tau)$ be a regular Cartan prolongation with
			an associated filtration by simple Cartan prolongations
			\[
				\pi^*I\subsetneq I_{(1)}\subsetneq\cdots\subsetneq I_{(L-1)}\subsetneq I_{(L)} = J.
			\]
			Let $I_k$ $(k = 0,1,\ldots,K)$ be the $k$-th extension of $I$ relative to $\pi$.
			Then we have 
				\[
					I_{(1)}\subset I_1.
				\]
		\end{lemma}
		\emph{Proof.} Since $\pi^*I\subset I_{(1)}$ represents a simple Cartan prolongation, 
				\[
					I_{(1)} = \mathcal{C}(\pi^*I)\cap I_{(1)} \subset \mathcal{C}(\pi^*I)\cap J = I_1,
				\]	
			as desired.	
			\qed	
		
		\begin{remark}
			It is not yet clear to us whether the conclusion in Lemma \ref{regFilVSrelExt} must hold for all $k$ in general.
		\end{remark}	
			
		\begin{Example}\label{singularExample}
			\emph{There exist singular Cartan prolongations.} 
			
			Consider two systems with the same independence condition $\ed t$:
			\begin{align*}
				I& = \lbb\ed x_1 - u_1\ed t,~ \ed x_2 - u_2\ed t\rbb,\\
				J& = \lbb I,~ \ed u_1 + f\ed u_2 - g \ed t,~ \ed f - h\ed t,~ \ed g - (f+h)\ed u_2\rbb.
			\end{align*}
			$I\subset J$ represents a Cartan prolongation. To see this, consider an integral curve 
			$\gamma: t\mapsto (\x, \u,t) =  (\x(t),  \x'(t), t)$ of $I$. A lifting of $\gamma$ to an
			integral curve of $J$ must satisfy:
			\begin{equation*}\left\{
				\begin{alignedat}{1}
				  u'_1 + f  u'_2   - g & = 0,\\
				  f' - h& = 0,\\
			 	 g' - (f+h) u'_2 & = 0.
				 \end{alignedat}\right.
			\end{equation*}
			From these equations, we obtain that
			\[
				(u'_1+ f u'_2)' = (f+ f')  u'_2.
			\]
			This implies that
			\begin{equation*}\left\{
				\begin{alignedat}{1}
				f &= \frac{u_1''}{u_2' - u_2''}~, \\
				g &= u_1'+ \frac{u_1''u_2'}{u_2' - u_2''}~, \\
				h &= \(\frac{u_1''}{u_2' - u_2''}\)'~,
				\end{alignedat}\right.
			\end{equation*}
			which is determined as long as $u_2' - u_2''\ne 0$ along $\gamma$.
			
		 	It is straightforward to compute that
				\[
					I_1 = \lbb I,~ \ed u_1+ f\ed u_2 - g\ed t\rbb.
				\]	
			$I\subset I_1$ is not a Cartan prolongation, since
			\[
				\rank(\mathcal{C}(I_1)) - \rank(\mathcal{C}(I)) = 7-5 = 2>1 = \rank(I_1) - \rank(I),
			\]
			violating Corollary \ref{simpleCPrankCond}. 
(The observation that $I \subset I_1$ is not a Cartan prolongation can also be seen more directly from the fact that a generic integral curve of $I$ does not have a unique lift to an integral curve of $I_1$.)
	Thus, the Cartan prolongation represented by
			$I\subset J$ is not $\mathcal{C}$-regular.
			
			Furthermore, if $I\subset J$ were regular with an associated filtration $I_{(\ell)}$ by  
			simple Cartan prolongations, then by Lemma \ref{regFilVSrelExt}, 
we would have			
$I_{(1)}\subset I_1$,
and
			 therefore $I_{(1)} = I_1$ (since $\rank(I_1/I) = 1$). However, this is impossible, since $I\subset I_1$
			does not represent a Cartan prolongation. 
			
					\end{Example}	
			
	\begin{Example}\label{regCregExample} \emph{There exist Cartan prolongations that are 
	regular but not $\mathcal{C}$-regular.}
	
	For $n\ge 3$, consider the following list of $1$-forms expressed in the 
	coordinates $(x_i, u_\alpha, v_\rho, w, t)$.
	\begin{equation*}\left\{
		\begin{alignedat}{1}
		\theta^i& = \ed x_i - u_i \ed t, \qquad (i = 1,\ldots,n)\\
		\eta^1& = \ed u_1 - v_1 \ed t,\\
		\eta^2& = \ed u_2 - v_2\ed u_3 - \cdots - v_{n-1}\ed u_n - w \ed t,\\
		\xi^1& = \ed v_1 - v_2 \ed t,\\
		&\quad \vdots\\
		\xi^{n-2}& = \ed v_{n-2} - v_{n-1}\ed t,\\
		\xi^{n-1}& = \ed v_{n-1} - v_n\ed t.
		\end{alignedat}\right.
	\end{equation*}
	Let
	\[
		I = \lbb\theta^1,\ldots,\theta^n\rbb, \quad J = \lbb I, \eta^1, \eta^2, \xi^1,\ldots,\xi^{n-1}\rbb.
	\]
	It is easy to see that $I\subset J$ represents a Cartan prolongation.
	Indeed, an integral curve $\gamma: t\mapsto (\x,\u,t)= (\x(t), \x'(t), t)$ of $I$ has its lifting to an integral curve
	of $J$ uniquely determined by the equations:
	\begin{equation*}\left\{
		\begin{alignedat}{1}
			v_1 &= u_1', \\
			 v_2 &= u_1'', \\
			 &\vdots\\
			  v_{n}& =  u_1^{(n)}, \\
			   w &= u_2'- v_2u_3'-\cdots- v_{n-1}u_n'. 
		\end{alignedat}\right.
	\end{equation*}
	Now $I_1 = \lbb I, \eta^1, \eta^2\rbb$, $I_2 = J$. Since $n\ge 3$, 
		\[
			\rank(\mathcal{C}(I_1)) - \rank(\mathcal{C}(I)) = n > 2 = \rank(I_1) - \rank(I), 
		\]
	and $I\subset I_1$ is not a Cartan prolongation.
	
	On the other hand, $I\subset J$ is regular, since we can take
	\begin{align*}
		I_{(0)}& = I,\\
		I_{(1)}& = \lbb I, \eta^1\rbb,\\
		I_{(2)}& = \lbb I, \eta^1, \xi^1\rbb,\\
			&\quad\vdots\\
		I_{(n)}& = \lbb I, \eta^1, \xi^1,\xi^2,\ldots, \xi^{n-1}\rbb,\\
		I_{(n+1)}& = \lbb I, \eta^1, \xi^1,\xi^2,\ldots, \xi^{n-1}, \eta^2\rbb = J,
	\end{align*}
	and $I_{(\ell)}$ $(\ell = 0,\ldots,n+1)$ provides a filtration of $J$ by simple Cartan prolongations.
	
	\end{Example}		
			
	The requirement ``$n\ge 3$" in Example \ref{regCregExample} is no coincidence, as the following
	theorem shows.
	
	\begin{theorem}\label{corank3Thm}
		Let $(M,I;\tau)$ be a system of corank $3$ (i.e., $\rank(\mathcal{C}(I)/I) = 3$). A 
		Cartan prolongation $\pi:(N,J;\sigma)\rightarrow (M,I;\tau)$ is regular if and only if
		it is $\mathcal{C}$-regular.
	\end{theorem}
	\emph{Proof.} $(\Leftarrow)$ is trivial; we now prove $(\Rightarrow)$ by induction on the order of $\pi$.
			If $\pi$ has order 1, then it is \emph{simple} by Lemma~\ref{rank1CPsimple}, and hence it is $\mathcal{C}$-regular.
			From now on, suppose that the theorem holds for Cartan prolongations of order less than $r$, and suppose that
			$\pi$ has order $r$.
			
			Let
			\[
				\pi^*I = I_0 \subsetneq I_1\subsetneq \cdots\subsetneq I_{K-1}\subsetneq I_K \subseteq J
			\]
		be the canonical filtration of $J$ by the relative extensions of $I$. Moreover, by assumption
		there exists a filtration
			\[
				\pi^*I = I_{(0)}\subsetneq I_{(1)}\subsetneq\cdots\subsetneq I_{(L-1)}\subsetneq I_{(L)} = J
			\]
		in which $I_{(\ell)}\subsetneq I_{(\ell+1)}$ are all simple Cartan prolongations.
		
		Suppose that $\rank (I_1/I) = 1$; then by Lemma~\ref{regFilVSrelExt}, $I_{(1)} = I_1$. Thus, 
		$I_1 = I_{(1)}\subset J$ represents a regular Cartan prolongation of order $r - 1$, which, by the inductive hypothesis,
		is $\mathcal{C}$-regular. It immediately follows that $\pi$ is also $\mathcal{C}$-regular.
		
		Thus, it suffices to consider the case when $\rank(I_1/I) > 1$. Since $(M,I;\tau)$ is a system of corank 3, 
		 $\rank(I_1/I)\le \rank(\mathcal{C}(I)) - \rank(I) = 3$. However, if $\rank(I_1/I) = 3$, then 
		$I_1$, and hence $J$, would contain the independence condition, violating the definition of a system. Therefore,
		$\rank(I_1/I) = 2$, which we assume from now on.
		
		Thus, by Lemma~\ref{regFilVSrelExt}, $\rank(I_{(1)}/I)$ is either $1$ or $2$.
		If $\rank(I_{(1)}/I) = 2$, then $I_{(1)} = I_1$. As in the case of $\rank(I_1/I) = 1$, one easily argues by induction
		that $\pi$ is $\mathcal{C}$-regular, so it remains to consider the case when $\rank(I_{(1)}/I) = 1$, which we now assume.
		
		Let $\eta$ be a nontrivial representative of $I_1/I_{(1)}$. In fact, suppose that, 
		in some coordinates, $I$ is spanned by
			\[
				\ed x_i - \sum_{j=1}^2 A_{ij}(\x, u_1, u_2, t)\ed u_j - B_i(\x, u_1, u_2, t)\ed t, \quad (i =  1,\ldots,n)
			\]		
		and that $I_{(1)}$ is spanned by $I$ and
 			\[
				\ed u_1 - f(\x,u_1,u_2, w, t)\ed u_2 - g(\x,u_1,u_2, w, t) \ed t,
			\]	
where $w$ is a fiber coordinate for the projection $M_{(1)} \to M$.
		We can choose 
			\[
				\eta = \ed u_2 - \lambda \ed t,
			\]
		 where $\lambda$ is independent of $\x, u_1, u_2, w,t$.
		 
		 Let $\mu$ be the smallest integer such
		 that $\eta$ is a section of $I_{(\mu)}$.
		 
		 Now consider the diagram below.
		 {\begin{figure*}[h!]
					\begin{center}
					\begin{tikzcd}[column sep = 10, row sep = 10]
						& \vdots \arrow[d, symbol = \supset]	\\
						& I_{(\mu)}\arrow[d, symbol = \supset]	 \arrow[ld, symbol = \supset]\\	
	\lbb I_{(\mu- 1)}, \eta\rbb\arrow[r, symbol = \supset]
		\arrow[d, symbol = \supset]	& I_{(\mu - 1)}\arrow[d, symbol = \supset]	\\	
			\vdots	\arrow[d, symbol = \supset]		&\vdots	\arrow[d, symbol = \supset]\\	
		\lbb I_{(1)}, \eta\rbb \arrow[r, symbol = \supset]
		\arrow[d, symbol = \supset]				&I_{(1)}\arrow[d, symbol = \supset]	\\	
				I_{(0)}\arrow[r, symbol = \eqdef]		&I							
 					\end{tikzcd}
					\end{center}
				\end{figure*}}	
		
		In this diagram, $I_{(0)}\subset \lbb I_{(1)}, \eta\rbb  = I_1$ represents a total prolongation. 
		Moreover, the $1$-form $\ed u_2$ must be independent of $I_{(k)}$ and $\ed t$ for $k = 1,\ldots, \mu - 1$;
		otherwise,
		either the minimality of $\mu$ would be violated, or $\ed t$ would be a section of $I_{(\mu)}$, 
		which is impossible.
		Thus, each horizontal inclusion in the diagram induces a prolongation by differentiation.
		This implies that
		$\lbb I_{(k)},\eta\rbb$ $(k = 1,\ldots, \mu - 1)$ are systems.
		
		Furthermore, for each $k \in \{1,\ldots,\mu - 2\}$, 
			\[
				\lbb I_{(k)}, \eta\rbb \subset \lbb I_{(k+1)}, \eta\rbb
			\]		
		represents a simple Cartan prolongation.  To see this, first note that the
		 underlying manifold of $\lbb I_{(k)},\eta\rbb$ is just $M_{(k)}\times \R$ with $\lambda$
		 as the coordinate on the $\R$-factor. Submersion at the manifold level for
		 consecutive $k$ follows. 
		 Because $\ed u_2$ is a section of 
		$\mathcal{C}(I_{(k)})$ and $\ed \lambda$ is not, we have
			\[
				\mathcal{C}(\lbb I_{(k)},\eta\rbb) \cap \lbb I_{(k+1)},\eta\rbb
					 = \lbb I_{(k+1)},\eta\rbb \quad (k = 1,\ldots, \mu - 2).
			\]
		Existence and uniqueness of lifting integral curves are evident.
				
		Now, regarding  $\rank(I_{(\mu)}/I_{(\mu - 1)})$, there are two possibilities---it is either $1$ or $2$. 
		(This is because $I_{(\mu - 1)}\subset I_{(\mu)}$ represents a simple Cartan prolongation of a system
		of corank $3$.) 
		
		\begin{enumerate}[(1)]
			\item{If $\rank(I_{(\mu)}/I_{(\mu - 1)}) = 1$, then $I_{(\mu)} = \lbb I_{(\mu -  1)}, \eta\rbb$;}	
			\item{If $\rank(I_{(\mu)}/I_{(\mu - 1)}) = 2$, that is, $I_{(\mu - 1)}\subset I_{(\mu)}$ represents a total prolongation. 
			Since we have argued that the inclusion $I_{(\mu - 1)} \subset \lbb I_{(\mu - 1)}, \eta\rbb$ is a prolongation by 
			differentiation, Theorem \ref{simpleCPThm} implies that 
						$\lbb I_{(\mu - 1)}, \eta\rbb\subset I_{(\mu)}$
							 represents a simple Cartan prolongation of order $1$.}
		\end{enumerate}
		In summary, we have the filtration
		\[
			\pi^*I = I_{(0)}\subsetneq \lbb I_{(1)}, \eta\rbb \subsetneq\cdots\subsetneq \lbb I_{(\mu- 1)}, \eta\rbb~\substack{\displaystyle\subsetneq\\ {\rm or} \\=} ~ I_{(\mu)}\subsetneq\cdots\subsetneq I_{(L)} = J,
		\]
		where each inclusion represents a simple Cartan prolongation.
		In other words, there exists a new filtration $I_{[k]}$ of $J$ by simple Cartan 
		prolongations that satisfies
		$\rank(I_{[1]}/I) = 2$, a case already treated.
		
		This completes the proof.\qed

		\begin{Prop}\label{kMinProp} 
		Let  $(M,I; \tau)$ be a system of corank 3, and let $\pi:(N,J;\sigma)\rightarrow (M,I;\tau)$ be 
		a $\mathcal{C}$-regular Cartan
		prolongation, with $I_k$ being the canonical relative extensions of $I$.
		Let $k$ denote the smallest integer such that $\rank(I_{k+1}/I_k)>1$. Then either $k$ does not exist, or 
		$k = 0$ and $\rank(I_1/I) = 2$.
		\end{Prop}
		\begin{proof}
		Suppose for the sake of contradiction that $k>0$.  Then
		$\rank(I_{\ell+1}/I_\ell) = 1$ for $\ell \in \{0,1,\ldots,k-1\}$. By the proof of Lemma \ref{regFilVSrelExt}, we have 
			\[
				I_{(\ell)} = I_\ell, \quad \ell = 0,1,\ldots,k.
			\]
		In other words, $I_1,\ldots,I_k$ are successive simple Cartan prolongations starting at $I$.  
		By Corollary \ref{simpleCPrankCond}, $(M_k,I_k;\sigma_k)$ is also a system of corank $3$. 
		It follows that 
			\[
				\rank(I_{k+1}/I_k)\le 2.
			\]
		By the characterization of $k$, we must have $\rank(I_{k+1}/I_k) = 2$.
		It follows that $I_{k+1}$ is equivalent to $\pr^{(1)}I_k$.
		
		As we will explain below, repeated application of Theorem \ref{simpleCPThm} yields the following diagram,
				{\begin{figure*}[h!]
					\begin{center}
					\begin{tikzcd}[column sep = 10, row sep = 5]
					& I_{k+1}\arrow[r, symbol = \cong]\arrow[dl, dashed, swap, "\rho"]\arrow[dd, "\beta"]& \pr^{(1)}I_k\\
		\pr^{(1)}I_{k-1}\arrow[dd,swap, "\phi"]\arrow[dr, "\psi"] && \\
				& I_k \arrow[dl, "\alpha"]	&\\
			 	   	I_{k-1}  	&&									
 					\end{tikzcd}
					\end{center}
				\end{figure*}}	
		where all arrows are understood as submersions at the manifold level.
		To be clear, in the diagram, the map
		$\alpha:M_k \to M_{k-1}$ is induced from the simple Cartan prolongation $I_{k-1}\subset I_k$, 
		and the map $\beta:M_{k+1} \to M_{k}$ is induced from the inclusion $I_k\subset I_{k+1}$, which represents a total 
		prolongation. Moreover, let $q_\ell: N\rightarrow M_\ell$ denote the quotient maps.
		
		Applied to $\alpha$, Theorem \ref{simpleCPThm} implies the existence of 
		 $\psi$, which is also a simple Cartan prolongation. Applied to $\psi$, the same theorem
		 implies the existence of $\rho$, again a simple Cartan prolongation. 

		 Because $\mathcal{C}$ commutes with pull-back, we have 
		 \[\pr^{(1)}I_{k-1}\subset \mathcal{C}(\phi^*I_{k-1}) = \phi^*\mathcal{C}(I_{k-1}).\]
		 Consequently, 
		  \[ 
		  	\begin{alignedat}{1}
				(\rho\circ q_{k+1})^*(\pr^{(1)}I_{k - 1}) 
				&\subset 
			 				(\rho\circ q_{k+1})^*(\phi^*\mathcal{C}(I_{k-1})) \\
					  &=(\alpha\circ\beta\circ q_{k+1})^*\mathcal{C}(I_{k-1}) 
			  		= q_{k-1}^* \mathcal{C}(I_{k-1}).
			  \end{alignedat}
		  \]
		  where the LHS is a subbundle of $J$, and the RHS is just $\mathcal{C}(I_{k-1})$, for $I_{k-1}\subset J$ is a subbundle.
		  To summarize, we have 
		 \[ (\rho\circ q_{k+1})^*(\pr^{(1)}I_{k - 1}) \subset \mathcal{C}(I_{k-1})\cap J = I_k, \]
		  This implies that the rank of $\pr^{(1)}I_{k-1}$ is at most the rank of $I_k$.
		   But this is impossible, since $\phi$ has order $2$, while $\alpha$ has order $1$.
			 
		 Therefore, if $k$ exists, then it must be zero, and $\rank(I_1/I) = 2$.
		 \end{proof}

		\begin{corollary}\label{totalFollowedbypartialCor}
			Let $(M,I;\tau)$ be a system of corank $3$. Suppose that $\pi: (N,J;\sigma)\rightarrow (M,I;\tau)$
			is a regular Cartan prolongation with the associated 
			 relative extensions $I_k$ of $I$. The chain of inclusions
				\begin{equation}\label{CartanFiltInCor}
					\pi^*I = I_0 \subsetneq I_1\subsetneq \cdots\subsetneq I_{K-1}\subsetneq I_K = J
				\end{equation}
			must represent a number (possibly zero) of successive total prolongations starting from $I$ followed by
			a number of successive simple Cartan prolongations of order $1$ that terminate at $J$.
		\end{corollary}
		
		\begin{proof}
			By Theorem~\ref{corank3Thm}, each inclusion in \eqref{CartanFiltInCor} is a simple Cartan prolongation.
			Furthermore, if $\rank(I_{k+1}/I_k)  \ge 2$ and $\rank(I_{\ell+1}/I_{\ell}) = 1$ for some $\ell < k$, 
			then the $\mathcal{C}$-regular Cartan prolongation $I_{\ell}\subset J$ would violate 
			Proposition~\ref{kMinProp}. 
		\end{proof}
		
		A particular case covered by Theorem \ref{corank3Thm} is when $\pi$ is obtained by 
		successive prolongations by differentiation: $\pi^*I\subset I_{(1)}\subset\cdots\subset I_{(L)} = J$. 
		When $\rank(I_{(1)}/I) = 2$, $I_{(1)} = I_1$ and $I_1\subset J$ represents a prolongation by successive differentiation (i.e., a prolongation obtained by performing a succession of prolongations by differentiation)
		whose order is less than that of $I\subset J$. When $\rank(I_{(1)}/I) = 1$, choose $\eta$ and determine the
		integer $\mu$ as in the proof of Theorem~\ref{corank3Thm}.
		Since, by assumption, each prolongation $I_{(k)}\subset I_{(k+1)}$ is obtained by differentiation, the same is true
		for $\lbb I_{(k)}, \eta\rbb\subset \lbb I_{(k+1)}, \eta\rbb$, $k\in\{1,\ldots,\mu-2\}$. 
		Now, if $\rank(I_{(\mu)}/I_{(\mu - 1)}) = 1$, then $I_{(\mu)} = \lbb I_{(\mu -  1)}, \eta\rbb$;
		if $\rank(I_{(\mu)}/I_{(\mu - 1)}) = 2$, then Theorem \ref{corank3Thm} and Remark \ref{simpleCPrmk}
		imply that 
						$\lbb I_{(\mu - 1)}, \eta\rbb\subset I_{(\mu)}$
		is a prolongation by differentiation.
		As a consequence,	
		\[
			\pi^*I = I_{(0)}\subsetneq \lbb I_{(1)}, \eta\rbb \subsetneq\cdots\subsetneq \lbb I_{(\mu- 1)}, \eta\rbb~\substack{\displaystyle\subsetneq\\ {\rm or} \\=} ~ I_{(\mu)}\subsetneq\cdots\subsetneq I_{(L)} = J
		\]
		is a filtration in which each inclusion represents a prolongation by differentiation, satisfying $\rank(\lbb I_{(1)}, \eta\rbb/I) = 2$.
		
		This argument justifies the following theorem.
		
		\begin{theorem}\label{corank3succdiffThm}
			Let $(M,I;\tau)$ be a system of corank $3$. If $\pi:(N,J;\sigma)\rightarrow (M,I;\tau)$ is 
			a prolongation by successive differentiation, then the associated  relative extensions $I_k$ of $I$
			satisfy: each $I_{k}\subset I_{k+1}$ represents a prolongation by differentiation. 
		\end{theorem}	 
		
		{\remark The conclusion in Theorem \ref{corank3succdiffThm} holds for general corank when all prolongations by differentiation
			are based on a fixed set of coordinates and the $t$-derivatives of the control variables
			therein. However, when the corank of $(M,I;\tau)$ is greater than $3$, the conclusion does not necessarily hold  if 
			one allows prolongation by differentiation to follow changes of coordinates. For example, consider
			\[
				I = \lbb \ed x_i - u_i \ed t\rbb_{i = 1}^3.
			\]
		The following is a prolongation of $I$ by successive differentiation\footnote{In fact, this prolongation is obtained by 
		setting $n = 3$ in Example \ref{regCregExample}.}:
			\begin{align*}
				J = \lbb I, ~& \ed u_1  - \alpha \ed t,\\
					& \ed \alpha - \beta\ed t,\\
					&\ed \beta - \gamma\ed t,\\
					& \ed(u_2 - \beta u_3) - w\ed t\rbb.	
			\end{align*}
		Direct calculation yields
			\[
				I_1 = \lbb I, \ed u_1 - \alpha\ed t, ~\ed u_2 - \beta \ed u_3 - \hat w\ed t\rbb,
			\]
		where $\hat w = \gamma u_3 +w$	. It turns out that $I_1$ is not even a Cartan prolongation of $I$, since, 
		otherwise, it must be a simple Cartan prolongation, but this is impossible by Corollary \ref{simpleCPrankCond}.

		}

	\subsection{Alternative Descriptions of $\ut$-Absolute Equivalence}  One can modify 
	the definition of $\uptau$-absolute equivalence by requiring the 
		Cartan prolongations to be regular or even to be
		ones obtained by successive differentiations. 
		The result is two new
		 equivalence relations among systems. It turns out that these new equivalence relations are 
		 no more restrictive
		than $\ut$-absolute equivalence, as we will demonstrate in this section.

	\begin{Def}
		Two systems $(M,I;\tau)$ and $(\bar M,\bar I;\bar\tau)$ 
		are said to be \emph{$\mathcal{R}$-related} (resp. $\mathcal{D}$-related) if there exists a system $(N,J;\sigma)$
	and submersions $\pi: N\rightarrow M$ and $\bar \pi:N\rightarrow \bar M$ that make $(N,J;\sigma)$
	a regular Cartan prolongation (resp., a prolongation by successive differentiation\footnote{Here we do not require 
	each prolongation by differentiation to be constructed from a fixed set of coordinates.}) of both $(M,I;\tau)$ 
	and $(\bar M,\bar I;\bar\tau)$.
	\end{Def}
	
	\begin{Prop}\label{DRAbsProp}
		Being $\mathcal{R}$-related (resp., $\mathcal{D}$-related) is an equivalence relation among systems.
	\end{Prop}
	\emph{Proof.} 
		Reflexivity and symmetry are trivial;  it suffices to prove transitivity.
		Suppose that $(M,I,\tau)$ and $(\bar M,\bar I;\bar\tau)$ are $\mathcal{R}$-related (resp.,  $\mathcal{D}$-related), 
		and suppose the same for $(\bar M,\bar I;\bar\tau)$  and
			$(\hat M,\hat I;\hat \tau)$. In the diagram below, by the assumption, $\pi,\bar \pi, \varpi,\hat\varpi$
			all represent regular Cartan prolongations (resp., prolongations by successive differentiation). 
			The rest of the diagram is constructed
			using the assumption that $\bar\pi$ and $\varpi$ are Cartan prolongations and Theorem \ref{SluisThm} (assuming $\bar L\ge L$); 
			in particular, $\pi_{\bar L,L}, \phi,\psi, \pi_{L,0}$ represent prolongations by successive differentiation.
				{\begin{figure*}[h!]
					\begin{center}
					\begin{tikzcd}[column sep = 15, row sep =6 ]
						&&\pr^{(\bar L)}I \arrow[d, swap, "\pi_{\bar L,L}"]\arrow[rdd, "\psi"]&&\\[0.8em]
						&& \pr^{(L)}\bar I\arrow[ld,swap, "\phi"]\arrow[dd,"\pi_{L,0}"]&&\\
						&J\arrow[ld,swap, "\pi"]\arrow[rd,swap,"\bar\pi"] && \bar J\arrow[ld,"\varpi"]\arrow[rd, "\hat\varpi"]&\\
						I && \bar I &&\hat I
 					\end{tikzcd}
					\end{center}
				\end{figure*}}
			The pair of maps $\pi\circ\phi\circ\pi_{\bar L, L}$ and $\hat\varpi\circ\psi$ thus
			both represent regular Cartan prolongations (resp., prolongations by successive differentiation). This completes the proof.\qed

	\begin{theorem}\label{RrelAbsEq}
		For systems, $\mathcal{D}$-related $\Leftrightarrow$ $\mathcal{R}$-related $\Leftrightarrow$  $\ut$-Absolutely equivalent.
	\end{theorem}
	\emph{Proof.} $(\Rightarrow)$ are obvious. `$\mathcal{D}$-related $\Leftarrow$ $\ut$-Absolutely equivalent' is a consequence of
			Sluis's extension theorem. In fact, suppose that $(N,J;\sigma)$ is a Cartan prolongation
			of both $(M,I;\tau)$ and $(\bar M, \bar I;\bar \tau)$ with the submersions $\pi$ and $\bar\pi$, respectively.
			Applying the diagram in
			Figure~\ref{Fig_SluisExt} to $\pi$,
			we note that all horizontal arrows (i.e., $\phi_k, \pi_k$, $k = 1,\ldots,K$) in that diagram
			represent prolongations by differentiation, and it follows that
			$(N,J;\sigma)$ and $(M,I;\tau)$ are $\mathcal{D}$-related. One similarly argues that
			$(N,J;\sigma)$ and $(\bar M,\bar I;\bar\tau)$ are also $\mathcal{D}$-related. Since being 
			$\mathcal{D}$-related is an equivalence relation (Proposition~\ref{DRAbsProp}), 
			$(M,I;\tau)$ and $(\bar M,\bar I;\bar \tau)$ 
			are $\mathcal{D}$-related.
			\qed\\

	{To end this section, we revisit the Cartan prolongations from Examples~\ref{singularExample}
	and \ref{regCregExample} and demonstrate explicitly that, in each case, the systems involved are $\mathcal{D}$-related.			
	\begin{Example}\label{example-4}
		Recall the singular Cartan prolongation in Example \ref{singularExample}:
		\begin{align*}
			I &= {\lbb \ed x_1 - u_1 \ed t, \ed x_2 - u_2 \ed t \rbb}\\
			&\subset {\lbb I, \ed u_1 + f\ed u_2 - g\ed t, \ed f - h\ed t, \ed g - (f+h)\ed u_2\rbb} = J.
		\end{align*}
		To establish a $\mathcal{D}$-relation between $I$ and $J$, it suffices to follow the construction leading to Theorem \ref{SluisThm}.
		Indeed, we prolong $J$ by differentiation:
		\begin{equation*}
			\begin{alignedat}{1}
			J_{(1)} &= \lbb J, \ed u_2 - \alpha \ed t\rbb, \\
			 J_{(2)} &= \lbb J_{(1)}, \ed \alpha - \beta \ed t\rbb,\\
			 J_{(3)} &= \lbb J_{(2)}, \ed \beta - \gamma \ed t\rbb.
			\end{alignedat}
		\end{equation*}
		And construct the $3$rd total prolongation of $I$:
		\begin{equation*}
			\begin{alignedat}{2}
			\pr^{(3)}I  = \lbb I,~& \ed u_1 - \lambda _1\ed t, ~&&\ed u_2- \lambda_2 \ed t,\\
							&\ed \lambda_1 - \mu_1 \ed t, &&\ed \lambda_2 - \mu_2\ed t, \\
							&\ed \mu_1 - \kappa_1\ed t,&& \ed \mu_2 - \kappa_2\ed t\rbb.
			\end{alignedat}
		\end{equation*}
		It is easy to determine a local isomorphism between $J_{(3)}$ and $\pr^{(3)} I$, restricted to a domain on which $\alpha\ne \beta$, as follows:
		\begin{equation*}\left\{
			\begin{alignedat}{1}
			\lambda_1& = g - f\alpha ,\\
			 \lambda_2 &= \alpha, \\
			 \mu_1 &= (\alpha-\beta)f,\\
			  \mu_2 &= \beta, \\
			  \kappa_1 &= (\alpha - \beta)h + (\beta - \gamma)f,\\
			   \kappa_2 &= \gamma,
			\end{alignedat}\right.
		\end{equation*}
	\end{Example}		
	
	\begin{Example}\label{example-5}
		For convenience, let $n = 3$ in Example \ref{regCregExample}. We have
			\begin{align*}
				I &= \lbb \ed x_i - u_i \ed t\rbb_{i = 1}^3\\
				& \subset \lbb I, ~\ed u_1- v_1 \ed t, ~\ed u_2 -v_2\ed u_3 - w\ed t, ~\ed v_1 - v_2\ed t
							,~\ed v_2 - v_3\ed t\rbb  = J.
			\end{align*}
		Following the construction in Theorem \ref{SluisThm}, we successively obtain
			\begin{align*}
				J_{(1)}&= \lbb J, \ed u_3 - \alpha \ed t\rbb,\\
				J_{(2)}& = \lbb J_{(1)}, \ed\alpha - \beta \ed t, \ed w - \gamma \ed t\rbb,\\
				J_{(3)}& = \lbb J_{(2)}, \ed (v_3\alpha + \gamma) - \eta\ed t, ~\ed\beta - \xi\ed t\rbb,
			\end{align*}
		and	
			\begin{equation*}
				\pr^{(3)}I  = \lbb I, \ed u_i - \lambda_i \ed t, ~\ed \lambda_i- \mu_i\ed t,~ \ed \mu_i - \kappa_i \ed t\rbb_{i = 1}^3.
			\end{equation*}
		A $\ut$-equivalence between $J_{(3)}$ and $\pr^{(3)}I$ can be established by the equations:
			\begin{equation*}\left\{
				\begin{alignedat}{1}
					\lambda_1 &= v_1, \\
					\lambda_2 &= v_2\alpha + w,\\
					\lambda_3 &= \alpha,\\
					 \mu_1 &= v_2,\\
					\mu_2 &= v_3\alpha+\gamma+v_2\beta, \\
					\mu_3 &= \beta, \\
					\kappa_1 &= v_3,\\
					 \kappa_2 &= \eta+v_3\beta + v_2\xi,\\
					\kappa_3 &= \xi.
				\end{alignedat}\right.
			\end{equation*}
		
		We point out that the prolongation by differentiation that generates $J_{(3)}$ from $J_{(2)}$ is not
		obtained by using the obvious coordinates in which
		$J_{(2)}$ is written; instead, it is obtained by first making a change of coordinates that turns $v_3\alpha+\gamma$ into a single variable.
	\end{Example}

%%Dynamic Feedback Linearization	
\section{$\ut$-Dynamic Linearization}\label{dynLinSec}

Given a control system, it is interesting to know whether we can transform it in a certain way 
into a (time-varying) linear system. When this is possible, such a transformation is often called 
a \emph{linearization} of the given system. 

The following notions of
\emph{linearization} are familiar in the literature. (See also \cite{de2018symmetry}.)

%%Comment
%{\color{red} 8. A citation here would probably be a good idea!\bf (Response: Citation added.)} {\color{blue} Great!}

\begin{enumerate}[\bf a.]
	\item{An autonomous control system $\dot \x = \f(\x,\u)$ is called \emph{static feedback linearizable} (SFL) if 
				there exists an invertible change of coordinates
				\begin{equation*}\left\{
					\begin{alignedat}{1}
					\y &= \phi(\x), \\
					 \v &= \psi(\x,\u),
					\end{alignedat}\right.
				\end{equation*}	
				that transforms the system into a linear system
				\[
					\dot \y = A\y + B\v,
				\]
				where $A,B$ are constant matrices.
				}
	\item{A time-varying control system $\dot \x = \f(t,\x,\u)$ is called \emph{extended static feedback linearizable} (ESFL) if
				there exists an ($t$-dependent) invertible change of coordinates
				\begin{equation}\label{tdeptrans}\left\{
					\begin{alignedat}{1}
					\y &= \phi(t,\x), \\
					\v &= \psi(t,\x,\u),
					\end{alignedat}\right.
				\end{equation}
				that transforms the system into a time-varying linear system
				\[
					\dot\y = A(t)\y+ B(t)\v.
				\]	
	}			
\end{enumerate}
 
In the generic case, one can find a
coordinate-independent criterion that works for both notions of linearizability above, as we now
explain.

\begin{Def}
	We say that a CTS $(M,I;\tau)$ is \emph{strongly linear} if 
	\begin{enumerate}[\bf i.]
		\item{the terminal derived system $I^{(\infty)} = 0$;}
		\item{each $\lbb I^{(k)},\tau\rbb$ is Frobenius.}
	\end{enumerate}
\end{Def}

{\remark \label{indCondRmk}
\begin{enumerate}[\bf A.]
\item{A corank-$p$ Pfaffian system $I$ corresponds to a distribution $\mathcal{D}$ on $M$.
	In the case of a CTS, Chow's theorem \cite{Chow39} implies that controllability corresponds to
	the bracket-generating property of $\mathcal{D}$, which is equivalent to
	 the condition $I^{(\infty)} = 0$.}
\item{Strong linearity is a property of a system $(M,I;\tau)$; in particular, it is sensitive to 
the independence condition $\tau$. Consider, for example, $(\R^4, I;\ed\alpha)$ with
	\[
		I = \lbb \ed f - g\ed \alpha,~ \ed g - h\ed \alpha\rbb
	\]
and $(\R^4, \bar I; \ed t)$ with
	\[
		\bar I = \lbb \ed x - \cos\theta \ed t, \ed y + \sin\theta \ed t\rbb.
	\]	
Via the diffeomorphism given by
	\begin{equation*}\left\{
	\begin{alignedat}{1}
		t& = f+h,\\
		x& =h \cos \alpha  + g\sin \alpha ,\\
		y& = -h\sin \alpha + g\cos \alpha,\\
		\theta & = \alpha,
	\end{alignedat}	\right.
	\end{equation*} 
$\bar I$ and $I$ correspond to each other. Therefore, the first derived systems
	\[
		I^{(1)} = \lbb \ed f - g \ed \alpha\rbb,\quad \bar I^{(1)} = \lbb \cos\theta \ed x - \sin\theta \ed y - \ed t\rbb
	\]
must also correspond under the diffeomorphism. 
However, $\lbb I^{(1)},\ed \alpha\rbb$ is integrable, while 
$\lbb \bar I^{(1)}, \ed t\rbb$ is not.  In other words, $(I;\ed \alpha)$ is strongly linear, while
$(\bar I;\ed t)$ is not.
}	
\end{enumerate}
}

\begin{theorem}	\cite{gardner1992gs, Sluis} \label{linFrob}
A controllable autonomous system $\dot\x =\f(\x,\u) $ is static feedback linearizable if and only if the corresponding
CTS is strongly linear.
\end{theorem}

\begin{theorem}{\rm(Cf. \cite[Theorem 3.11]{de2018symmetry})}
Let $\dot \x = \f(t,\x,\u)$ be a controllable system with $n$ states and $m$ inputs. 
Let $(M,I,\ed t)$ denote the corresponding 
CTS. The following are equivalent:
	\begin{enumerate}[\bf i.]
		\item{$(M,I,\ed t)$ is strongly linear;}
		\item{there exists a transformation \eqref{tdeptrans} that turns the system into 
					a time-varying linear system \[\dot \y = A(t)\y+B(t)\v;\]}
		\item{there exists a transformation \eqref{tdeptrans} that turns the system into
					a Brunovsk\'y normal form.}
	\end{enumerate}
\end{theorem}
			
\begin{proof}
	It suffices to prove \textbf{ii}$\Rightarrow$\textbf{i}$\Rightarrow$\textbf{iii}.
	
	(\textbf{ii}$\Rightarrow$\textbf{i}) Consider a time-varying CTS $(M,I;\ed t)$ where $I$
	is spanned by the $n$ $1$-forms\footnote{We adopt the convention of summing over repeated indices.}
	\[ 
		\theta^i = \ed x^i - \left(A^i_j(t) x^j + B^i_\alpha (t) u^\alpha\right)\ed t,
	\]
	where $i,j = 1,\ldots, n$ and $\alpha,\beta = 1,\ldots, m$. Let $n_1:= n - m$. Without loss of generality, assume
	that the $m$-by-$m$ minor $\det(B^{n_1+\alpha}_\beta)\ne 0$. Thus, by a change of coordinates of the form
	$u^\alpha \mapsto Q^\alpha_\beta(t)u^\beta$, we can arrange that 
	$B^{n_1+\alpha}_\beta = \delta^\alpha_\beta$.
	Using this, one easily finds that $I^{(1)}$ is spanned by
	\[
		\eta^\rho = \theta^\rho -  B^\rho_\alpha(t) \theta^{n_1 + \alpha}
		\qquad (\rho = 1,\ldots, n_1).
	\] 
	By introducing $\bar x^\rho  = x^\rho -  B^\rho_\alpha(t) x^{n_1+\alpha}$, we find that
	each $\eta^\rho$ is of the form 
	\[
		\eta^\rho = \ed\bar x^\rho - \left[A^\rho_\sigma(t) \bar x^\sigma + C^\rho_\alpha(t) x^{n_1+\alpha}\right] \ed t
	\]
	for some functions $C^\rho_\alpha(t)$. Thus, $I^{(1)}$ is also a time-varying linear system, 
	and $\lbb I^{(1)}, \ed t\rbb$ is Frobenius. This procedure continues. Controllability implies that $I^{(\infty)} = 0$.
	Therefore, $I$ is strongly linear.
	
	(\textbf{i}$\Rightarrow$\textbf{iii}) Let $(M,I;\tau)$ be a strongly linear CTS with $I^{(K-1)}\ne 0$ and $I^{(K)} = 0$.
	Suppose that $\tau = \ed t$, and let $s_1 = \rank(I^{(K-1)})$. By assumption, $\lbb I^{(K-1)}, \ed t\rbb$ is Frobenius.
	Since $\ed t$ does not belong to $I$, we have
	\[
		 I^{(K-1)} = \lbb \ed x^{K-1}_i - x_i^{K-2}\ed t\rbb_{i = 1}^{s_1}
	\]
	for some functions $x^{K-1}_1, \ldots, x^{K-1}_{s_1}, x^{K-2}_{1}, \ldots, x^{K-2}_{s_1}$; these $2s_1$ functions and $t$
	have linearly independent differentials, because $I^{(K)} = 0$.
	
	Now, suppose that for some $\ell \ge 0$ we have shown (as we just did for $\ell = 1$) that
	\begin{equation}\label{IKell}
		I^{(K - \ell - 1)} = \lbb I^{(K - \ell)}, \ed x_{i}^{K - \ell - 1} - x_i^{K - \ell - 2}\ed t\rbb_{i = 1}^{s_\ell},
	\end{equation}
	where $\ed x_i^{K-\ell - 1}$ and $\ed x_i^{K - \ell - 2}$ are independent of $\mathcal{C}(I^{(K-\ell)})$ and among themselves.
	Since $I^{(K - \ell - 1)}$ is the derived system of $I^{(K - \ell - 2)}$, we have
	\[
		\ed(\ed x_i^{K - \ell - 1} - x_i^{K - \ell - 2}\ed t) = \ed t\W \ed x_i^{K - \ell - 2} \equiv 0 \mod I^{(K - \ell - 2)}. 
	\]
	Combined with the assumption that $\lbb I^{(K - \ell - 2)}, \ed t\rbb$ is Frobenius, we see that there exist new functions
	$x_{1+s_\ell}^{K - \ell - 2}, \ldots, x_{s_{\ell +1}}^{K - \ell - 2}$ and $x_{s_1}^{K - \ell - 3}, \ldots, x_{s_{\ell +1}}^{K - \ell - 3}$
	such that
	\[
		I^{(K - \ell - 2)} = \lbb I^{(K - \ell - 1)}, \ed x_i^{K - \ell - 2} - x_i^{K - \ell - 3}\ed t \rbb_{i = 1}^{s_{\ell +1}}
	\]
	where $s_{\ell+1} = \rank(I^{(K - \ell - 2)}) - \rank(I^{(K - \ell - 1)}) \ge s_\ell$. Because no combination of the $\ed x_i^{K - \ell - 2} - x_i^{K - \ell - 3}\ed t$ occurs in the derived system $I^{(K - \ell - 1)}$, the differentials $\ed x_i^{K - \ell - 2}$ and $\ed x_i^{K - \ell - 3}$ ($i = 1,\ldots, s_{\ell+1}$) must be independent of $\mathcal{C}(I^{K - \ell - 1})$ and among themselves.
	
	Repeat this procedure until $\ell = K - 1$ in \eqref{IKell},  at which point one recognizes that $I$ is in a Brunovsk\'y normal form.
\end{proof}			
			
		\begin{Def}
			A CTS is called \emph{$\ut$-dynamically linearizable} if it is $\ut$-absolutely equivalent
			to a strongly linear CTS.
		\end{Def}	
			
		The following theorem and corollary, slightly modified from their original statements
		in \cite{Sluis}, reduces 
		the problem of dynamic feedback linearization to finding a particular 
		type of Cartan prolongation.
			
		\begin{theorem}\cite[p.87, Theorem 41]{Sluis}\label{totprolLinear}
			A CTS is strongly linear if and only if its total prolongation is strongly linear. 
		\end{theorem}	
		\emph{Proof.} Let $(M,I;\tau)$ be a CTS. It is easy to verify that 
			\[
				(\pr^{(1)}I)^{(1)} = I.
			\]	
			The conclusion follows by the definition of strong linearity.\qed
			
		\begin{corollary}\label{dynlinCorollary}
			A CTS is $\ut$-dynamically linearizable if and only if 
			there exists a prolongation by successive differentiation that results in 
			a strongly linear system.
		\end{corollary}	
		{\begin{figure*}[h!]
					\begin{center}
					\begin{tikzcd}[column sep = 13, row sep =6 ]
						&& \pr^{(L)}\bar I\arrow[ld,swap, "\phi"]\arrow[dd,"\pi_{L,0}"]\\
						&J\arrow[ld,swap, "\pi"]\arrow[rd,swap,"\bar\pi"] &\\
						I && \bar I 
 					\end{tikzcd}
					\end{center}
				\end{figure*}}
		\emph{Proof.} 
			For $(\Leftarrow)$, assume that $J$ is obtained from $I$ by successive prolongation by differentiation
			and that $J$ is strongly linear. This implies that $I$ and $J$ are absolutely equivalent, 
			and by definition $I$ is $\ut$-dynamically linearizable. 
			For $(\Rightarrow)$, suppose that $(M,I;\tau)$
			and $(\bar M,\bar I;\bar\tau)$ are $\uptau$-absolutely equivalent, where
			$(\bar M,\bar I;\bar\tau)$ is strongly linear. Theorem \ref{RrelAbsEq} implies that
			these two systems are also $\mathcal{D}$-related, which yields $\pi$ and $\bar\pi$
			(prolongations by successive differentiation)
			in the diagram above.
			Now, for $\bar\pi$, we apply Theorem~\ref{SluisThm}, which yields the prolongations 
			$\phi$ and $\pi_{L,0}$. All the arrows in the diagram represent prolongations by differentiation.
			Thus, the same is true for $\pi\circ\phi$. The system $\pr^{(L)}\bar I$ is strongly linear
			by Theorem \ref{totprolLinear}. This completes the proof.\qed

		\begin{lemma}\label{typen2dynLinLemma}
			A type $(n,2)$ CTS $(M,I;\tau)$ is $\ut$-dynamically linearizable
			if and only if it admits a $\mathcal{C}$-regular Cartan prolongation $(N,J;\sigma)$
			that satisfies
			\begin{enumerate}[\bf i.]
				\item{$(N,J;\sigma)$ is strongly linear;}
				\item{each relative extension $I_k$ is a CTS}.
			\end{enumerate}
		\end{lemma}
		\emph{Proof.} For $(\Leftarrow)$, since $J$ is a Cartan prolongation of $I$, 
		they are absolutely equivalent, and since
		$J$ is assumed to be strongly linear, $I$ is $\ut$-dynamically linearizable by definition.
		For $(\Rightarrow)$, Corollary \ref{dynlinCorollary} 
		and Theorem \ref{corank3succdiffThm} together imply that there exists a Cartan prolongation
		terminating at a strongly linear system
		with each associated $I_{k+1}$ being a prolongation by differentiation of $I_k$.
		Such a Cartan prolongation is automatically $\mathcal{C}$-regular and satisfies
		{\bf i} and {\bf ii}, since a prolongation by differentiation of a CTS results in a CTS.
		\qed
		 	
		\begin{lemma}\label{noTotLemma}
			Let $(M,I;\tau)$ be a type $(n,2)$ CTS. Suppose that 
			$\pi: (N,J;\sigma)\rightarrow (M,I;\tau)$ is a $\mathcal{C}$-regular
			Cartan prolongation (with extension length $K$) that satisfies:
				\begin{enumerate}[\bf i.]
					\item{$(N,J;\sigma)$ is strongly linear;}
					\item{each relative extension $I_k$ is a CTS;}
					\item{$\rank(I_1/I) = 2$;}
					\item{$\rank(I_{k+1}/I_{k}) = 1$ for all $1\le k\le K-1$.}
				\end{enumerate}
			Then $I\subset J^{(1)}$ represents
			a $\mathcal{C}$-regular Cartan prolongation satisfying conditions {\bf i} and {\bf ii}.
		\end{lemma}	
		\emph{Proof.} 	The assumptions imply that one can choose a coframing 
				\[(\theta^1,\ldots,\theta^n, \xi^1,\eta^1,\ldots,\eta^K,\omega^1,\omega^2, \tau)\]
				on $N$ such that (after dropping pull-back symbols):
			\begin{align*}
				I & = \lbb \theta^1,\ldots,\theta^n \rbb,\\
				I^{(1)}& = \lbb \theta^2,\ldots,\theta^{n-1}\rbb,\\
				I_k& = \lbb \theta^1, \ldots,\theta^n, \xi^1,\eta^1,\ldots,\eta^k\rbb, \quad k = 1,\ldots,K,\\
				\mathcal{C}(I_1) &= \lbb I_2,\omega^1, \tau\rbb.
			\end{align*}
			
		In particular, $I_K = J$, by the definition of $\mathcal{C}$-regularity.
			
	%%Comment
	%{\color{red} 9. Maybe comment explicitly that $I_K = J$?  Also, this is another place where I feel like a very simple example in coordinates would be illuminating! \bf (Response: Sentence added before. Example added after this proof.)} {\color{blue} Looks great!}
			
			Next, it is not difficult to see that such a coframing can be chosen to further satisfy
			the structure equations:
			\begin{equation}\label{CTSlinLemmaCong}\left\{
				\begin{alignedat}{2}
				\ed\theta^1& \equiv\tau\W\xi^1&&\mod I,\\
				\ed\theta^\alpha& \equiv 0&&\mod I,	\quad(\alpha = 2,\ldots,n-1)\\
				\ed\theta^n& \equiv \tau\W\eta^1&&\mod I,\\
				\ed\xi^1& \equiv \tau \W\omega^1&&\mod I_1,\\
				\ed\eta^k &\equiv \tau\W\eta^{k+1} &&\mod I_k,\quad  (k = 1,\ldots,K-1)\\
				\ed\eta^K&\equiv \tau\W\omega^2&&\mod I_K.
				\end{alignedat}\right.
			\end{equation}
			Note, in particular, that each $I_k$ being a CTS enforces that, in the 
			congruences above, the right-hand-sides
			are multiples of $\tau$.

			Note that
			\[
				J^{(1)} = \lbb I, \eta^1,\ldots,\eta^{K-1}\rbb.
			\]
			By the assumption {\bf i}, $\lbb J^{(1)}, \tau\rbb$ is Frobenius. It follows that there exist
			functions $A^k$ on $N$ such that 
			\[
				\ed\eta^k \equiv \tau\W\eta^{k+1} + A^k \tau\W\xi^1\mod I, \eta^1,\ldots,\eta^k
			\]
			for $k = 1,\ldots,K-1$. In fact, we can arrange all $A^k$ $(k = 1,\ldots,K-1)$ to be zero
			by adding an appropriate multiple of $\theta^1$ into each $\eta^k$.
			
			It now follows from the congruences \eqref{CTSlinLemmaCong}
			  that $I\subset J^{(1)}=: \bar J$ 
			represents a $\mathcal{C}$-regular Cartan
			prolongation of $I$ with the relative extensions
			\[
				\bar I_k = \lbb I, \eta^1,\ldots,\eta^k\rbb.
			\]		
			It is clear that this Cartan prolongation satisfies {\bf i} and {\bf ii}, and its
			extension length is $K-1$.
						\qed\\
				
		\begin{Example} Consider the type $(3,2)$ CTS $(M,I;\ed t)$ generated by
			\[
				\left\{
					\begin{alignedat}{1}
						\theta^1 & = \ed x_1 - (x_2+ uv)\ed t,\\
						\theta^2 & = \ed x_2 - (u+x_1v)\ed t,\\
						\theta^3& = \ed x_3 - v\ed t.
					\end{alignedat}
				\right.
			\]
			Let $J = \lbb I, \theta^4,\ldots, \theta^7\rbb$, where
			\[
				\left\{
					\begin{alignedat}{1}
						\theta^4 & = \ed u - u_1\ed t,\\
						\theta^5 & = \ed v - v_1\ed t,\\
						\theta^6& = \ed v_1 - v_2\ed t,\\
						\theta^7& = \ed v_2 - v_3\ed t.
					\end{alignedat}
				\right.
			\]
			
			It is observed in \cite[Example 48]{Sluis} that the inclusion $I\subset J$ induces a 
			Cartan prolongation, and $(J;\ed t)$ is strongly linear. It is easy to see that 
			\[
				I_1 = \lbb I, \theta^4,\theta^5\rbb, \quad I_2 = \lbb I_1, \theta^6\rbb,
				\quad I_3 = \lbb I_2, \theta^7\rbb = J.
			\]
			Lemma  \ref{noTotLemma} applies, 
			indicating that $I\subset \lbb I, \theta^5,\theta^6\rbb$ also represents a 
			 $\mathcal{C}$-regular
			Cartan prolongation with $(\lbb I, \theta^5,\theta^6\rbb; \ed t)$ being strongly linear.
			A direct verification of this fact will be left to the interested reader. We note that this linearization
			of $(M,I;\tau)$ is different from the one given by \cite[Example 48]{Sluis}; the latter is 
			obtained by differentiating $u$ instead of $v$.

		\end{Example}								
		
		The following theorem will serve as a basis for classifying $\ut$-dynamically linearizable
		type $(n,2)$ systems.				
						
		\begin{theorem}\label{typen2LinTheorem}
			A type $(n,2)$ CTS $(M,I;\tau)$ is $\ut$-dynamically linearizable if 
			and only if either it is already strongly linear, or it admits a $\mathcal{C}$-regular Cartan prolongation 
			$(N,J;\sigma)$ (with extension length $K$)
			satisfying:
				\begin{enumerate}[\bf i.]
					\item{$(N,J;\sigma)$ is strongly linear;}
					\item{each relative extension $I_k$ is a CTS;}
					\item{$\rank(I_{k+1}/I_k) =1$, for all $k = 0,\ldots,K-1$.}
				\end{enumerate}
		\end{theorem}	
		\emph{Proof.} By Lemma \ref{typen2dynLinLemma}, it suffices to justify {\bf iii}
		in $(\Rightarrow)$. Suppose that we already have a $\mathcal{C}$-regular 
		Cartan prolongation
		of $(M,I;\tau)$ satisfying {\bf i} and {\bf ii}. By Corollary \ref{totalFollowedbypartialCor}, 
		as $k$ increases from $0$, 
		the simple Cartan prolongations represented by $I_{k}\subsetneq I_{k+1}$
		are a number of (if any) successive total prolongations followed by
		Cartan prolongations of order 1. If there is any total prolongation in this list, we can apply
		Lemma \ref{noTotLemma} to find a $\mathcal{C}$-regular Cartan prolongation 
		of $(M,I;\tau)$ with a lower order and still satisfying {\bf i} and {\bf ii}. Continue until
		either $J = I$ or none of $I_k\subsetneq I_{k+1}$ represents a total prolongation; the former
		case implies that $I$ is strongly linear, by Theorem~\ref{totprolLinear} and Lemma~\ref{noTotLemma}, and the latter case
		implies \textbf{iii}.
		\qed
			
		\begin{Def}
			We say that a type $(n,2)$ CTS $(M,I;\tau)$ has \emph{class $K$} if 
			it is $\ut$-dynamically linearizable with $K$ being the minimal 
			integer such that there exists a $\mathcal{C}$-regular Cartan prolongation (with 
			extension length $K$)
			of $(M,I;\tau)$ satisfying the conditions {\bf i}-{\bf iii} in Theorem \ref{typen2LinTheorem}.
			If $(M,I;\tau)$ is not $\ut$-dynamically linearizable, we say that it
			has \emph{class $\infty$}.
		\end{Def}	
			
		\begin{theorem}\label{typen2StrEqn}
			\begin{enumerate}[\bf a.]
				\item{Suppose that $(M,I;\tau)$ is $\ut$-dynamically linearizable. If
			$\pi:(N,J;\sigma)\rightarrow (M,I;\tau)$ is a $\mathcal{C}$-regular Cartan
			prolongation (with extension length $K$) satisfying the
			conditions {\bf i}-{\bf iii} in Theorem \ref{typen2LinTheorem}, then on $N$ there exists
			a local coframing 
				\begin{equation}\label{0adaptedcoframing}
					(\theta^1,\ldots,\theta^n,\eta^1,\ldots,\eta^K,\omega^1,\omega^2, \sigma)
				\end{equation}
			satisfying
				\begin{align}
					I &= \lbb \theta^1,\ldots,\theta^n\rbb,\label{0adaptedCond1}\\
					I_k &= \lbb I, \eta^1,\ldots,\eta^k\rbb, \quad k = 1,\ldots,K,	\label{0adaptedCond2}
				\end{align}	
			and the structure equations:
				\begin{equation}\label{0adaptedstreq}\left\{
					\begin{alignedat}{2}
					\ed\theta^1& \equiv\sigma\W\omega^1&&\mod \theta^1,\ldots,\theta^n,\\
					\ed\theta^\alpha& \equiv 0&&\mod \theta^1,\ldots,\theta^n,	\quad(\alpha = 2,\ldots,n-1)\\
					\ed\theta^n& \equiv \sigma\W\eta^1&&\mod \theta^1,\ldots,\theta^n,\\
					\ed\eta^k &\equiv \sigma\W\eta^{k+1} &&\mod \theta^1,\ldots,\theta^n,\eta^1,\ldots,\eta^k,\quad  (k = 1,\ldots,K-1)\\
					\ed\eta^K&\equiv \sigma\W\omega^2&&\mod \theta^1,\ldots,\theta^n,\eta^1,\ldots,\eta^K.\\
					\end{alignedat}\right.
				\end{equation}	}
			\item{Conversely, if $(N,J;\sigma)$
			is a strongly linear system
			with a coframing \eqref{0adaptedcoframing} satisfying
				\[
					J = \lbb \theta^1,\ldots,\theta^n,\eta^1,\ldots,\eta^K\rbb	
				\]
			and the structure equations \eqref{0adaptedstreq}, then the class of the  system
				\[
					I:= \lbb \theta^1,\ldots,\theta^n\rbb
				\]
			(with independence condition induced by $\sigma$) is at most $K$.}
			\end{enumerate}	
		\end{theorem}	
		\emph{Proof.} {\bf a.} First, by the assumption, it is clear that there exists a coframing
			\eqref{0adaptedcoframing} satisfying \eqref{0adaptedCond1}, 
			\eqref{0adaptedCond2} and 
			\[
				\mathcal{C}(I) = \lbb I,\eta^1,\omega^1, \sigma\rbb.
			\] 
			Since $I$ is a CTS with corank $3$, for each $i = 1,\ldots,n$, 
			there must exist functions $A^i, B^i$ such that 
			\[
			 	\ed\theta^i \equiv \sigma\W(A^i\eta^1+B^i\omega^1)\mod I.
			\]
			 Moreover, the $n\times 2$ matrix $(A^i|B^i)$ must have rank $2$.
			It follows that one can make a linear transformation of the $\theta^i$ to arrange that
			\[
				A^n = B^1 = 1
			\]
			and all other $A^i, B^i$ are zero.
		
			Continuing, by the construction of $I_2$ and the expression for $\ed \theta^1$, we have
				\[
					\mathcal{C}(I_1) \subseteq \lbb I, \eta^2,\omega^1,\sigma\rbb.
				\]
			This
			inclusion must be an equality
			because $I\subset I_1$ represents a simple Cartan prolongation, which preserves 
			corank (Corollary \ref{simpleCPrankCond}). Since $I_1$, by assumption, is a
			CTS, it follows that 
				\[
					\ed\eta^1 \equiv \sigma\W(C^1\eta^2+ D^1\omega^1) \mod I_1
				\]
			for some functions $C^1,D^1$. We can always arrange that $D^1=0$ by adding
			a multiple of $\theta^1$ to $\eta^1$. On the other hand, $C^1$ is nonvanishing;
			hence,  by scaling $\eta^2$, we can arrange that $C^1 = 1$.
			
			We can continue with this type of argument and obtain the congruences
				\[
					\ed\eta^k \equiv \sigma\W\eta^{k+1} \mod I_k
				\]
			for $k = 2,\ldots,K-1$.
			
			Now, since $\mathcal{C}(I_K) = T^*N$ and $I_K = J$ is a CTS, it follows that, modulo $J$,
			$\ed\eta^K$ is congruent to a linear combination of 
			$\sigma\W\omega^i$ $(i = 1,2)$.
			We can add a multiple of $\theta^1$ to $\eta^K$ and then scale $\omega^2$ to arrange 
			that
			 	\[
					\ed\eta^K\equiv\sigma\W\omega^2\mod J.	
				\]	
			Thus we have obtained a desired coframing on $N$.\\
			
			{\bf b.} Assuming a coframing \eqref{0adaptedcoframing} on $N$ satisfying \eqref{0adaptedstreq},
				and letting 
					\[
						I = \lbb\theta^1,\ldots,\theta^n\rbb,
					\]
				it is easy to see that
				\[
					\mathcal{C}(I) = \lbb I, \eta^1, \omega^1,\sigma\rbb.
				\]
				Furthermore, $\lbb I,\sigma\rbb$ is clearly Frobenius.
				It follows that the system $I$ is a corank $3$ CTS.
				
				It is easy to see that $I\subset J$ represents a $\mathcal{C}$-regular 
				Cartan prolongation with
				\[
					I_k = \lbb I, \eta^1,\ldots,\eta^k\rbb, \quad k = 1,\ldots,K.
				\]
				By the assumption that $(N,J;\sigma)$ is strongly linear and the structure equations, the
				conditions {\bf i}-{\bf iii} in Theorem \ref{typen2LinTheorem}
				are satisfied. 
				
				This completes the proof. \qed\\

				Theorem \ref{typen2StrEqn} enables us to take the following approach towards
				 finding type $(n,2)$ CTS that are
				$\ut$-dynamically linearizable: \emph{For each $n, K$, first find 
				local coframings \eqref{0adaptedcoframing} adapted to some
				strongly linear system $(N,J;\sigma)$ 
				such that the structure equations \eqref{0adaptedstreq} hold; then we classify
				the CTS generated by $\theta^1,\ldots,\theta^n$ with
				independence conditions induced by $\sigma$.} 
								
				A concrete discussion of the classification problem, particularly in the case when 
				$n = 3$, will be addressed in a sequel to the current paper.

\section{Acknowledgement}						

We would like to thank Prof. Robert L. Bryant for generously offering his time 
for discussion.
This work also benefited from discussions with Prof. George R. Wilkens and Taylor J. Klotz.	

Additionally, we would like to thank the referee for their comments and questions, which have led to substantial improvements in our exposition.

\bibliographystyle{alpha}
\normalbaselines

\end{document}